\documentclass[11pt]{article}

\usepackage[small,bf]{caption}
\usepackage{amsmath, amsthm, amssymb}
\usepackage{algorithmic}
\usepackage{algorithm}
\usepackage[parfill]{parskip}

\usepackage{epsfig,fullpage,psfrag,url,verbatim}
\usepackage{graphics}
\usepackage{graphicx}
\usepackage{amsmath}
\usepackage{titlesec}
\usepackage{lipsum}

\usepackage{amsfonts}

\usepackage{color}

\newtheorem{thm}{Theorem}[section]

\newtheorem{prop}{Proposition}[section]

\newcommand{\inside}{\mathrm{int}}

 \medskip

\newcommand{\ua}{f_{\mathrm{slb}}}
\newcommand{\optim}{\mathrm{Opt}}
\newcommand{\slice}{\mathrm{Slice}}
\newcommand{\dist}{\mathrm{Dist}}
\newcommand{\diam}{\mathrm{Diam}}
\newcommand{\by}{\mathbf{y}}

\newcommand{\bX}{\mathbf{X}}

\DeclareMathOperator*{\argmin}{arg\,min}

\allowdisplaybreaks

\titlespacing*{\section}
{0pt}{1.6ex }{1.6ex}
\titlespacing*{\subsection}
{0pt}{1ex }{1ex}

\begin{document}

\title{New Computational Guarantees for Solving Convex Optimization Problems with First Order Methods, via a Function Growth Condition Measure}
\author{Robert M. Freund\thanks{MIT Sloan School of Management, 77 Massachusetts Avenue, Cambridge, MA   02139
({mailto:  rfreund@mit.edu}).  This author's research is supported by AFOSR Grant No. FA9550-15-1-0276 and  the MIT-Belgium Universit\'{e} Catholique de Louvain Fund.}
\and Haihao Lu\thanks{MIT Department of Mathematics, 77 Massachusetts Avenue, Cambridge, MA   02139
({mailto:  haihao@mit.edu}).}}
\date{Revised October, 2016} 

\maketitle

\begin{abstract} Motivated by recent work of Renegar \cite{renegar2015}, we present new computational methods and associated computational guarantees for solving convex optimization problems using first-order methods.  Our problem of interest is the general convex optimization problem $f^* = \min_{x \in Q} f(x)$, where we presume knowledge of a strict lower bound $\ua < f^*$.  [Indeed, $\ua$ is naturally known when optimizing many loss functions in statistics and machine learning (least-squares, logistic loss, exponential loss, total variation loss, etc.) as well as in Renegar's transformed version of the standard conic optimization problem \cite{renegar2015}; in all these cases one has $\ua = 0 < f^*$.] We introduce a new functional measure called the growth constant $G$ for $f(\cdot)$, that measures how quickly the level sets of $f(\cdot)$ grow relative to the function value, and that plays a fundamental role in the complexity analysis.  When $f(\cdot)$ is non-smooth, we present new computational guarantees for the Subgradient Descent Method and for smoothing methods, that can improve existing computational guarantees in several ways, most notably when the initial iterate $x^0$ is far from the optimal solution set.  When $f(\cdot)$ is smooth, we present a scheme for periodically restarting the Accelerated Gradient Method that can also improve existing computational guarantees when $x^0$ is far from the optimal solution set, and in the presence of added structure we present a scheme using parametrically increased smoothing that further improves the associated computational guarantees.\end{abstract}

\section{Problem Statement and Overview of Results}\label{intro}

\subsection{Problem Statement, Strict Lower Bound, and Function Growth Constant}

Motivated by recent work of Renegar \cite{renegar2015}, we present new computational methods and associated computational guarantees for solving convex optimization problems using first-order methods.  Our problem of interest is the following optimization problem:
\begin{equation}\label{poi1}
\begin{array}{lrlr} P:  & f^* := \ \  \mbox{minimum}_x & f(x)\\
& \mbox{ s.t. } & x \in Q \ ,
\end{array}
\end{equation}
where $Q \subseteq \mathbb{R}^n$ is a closed convex set and $f(\cdot) : Q \to \mathbb{R}$ is a convex function.  Let the set of optimal solutions of \eqref{poi1} be denoted as $\optim := \{ x \in Q : f(x) =f^* \}$.  For $x \in Q$, let $ \dist (x,\optim )$ denote the distance from $x$ to the set of optimal solutions, namely $ \dist (x,\optim ) := \min_y\{\|y-x\| : y \in \optim\}$.

{\bf Strict Lower Bound $\ua$ and Function Growth Constant $G$.}  Let $\ua$ be a known and given strict lower bound on the optimal value $f^*$ of \eqref{poi1}, namely  $\ua < f^*$.  Such a known strict lower bound arises naturally when optimizing many loss functions in statistics and machine learning (least-squares loss, logistic loss, exponential loss, total variation loss, etc.) perhaps with the addition of a regularization term; in all these cases $\ua = 0 < f^*$.  A known strict lower bound also arises in Renegar's transformed version of the standard conic optimization problem \cite{renegar2015}.

 Let $\varepsilon' >0$ be given.  Given the knowledge of the strict lower bound $\ua$, it is natural to work with the notion of a relative measure of optimality.  Let us define an $\varepsilon'$-relative solution of \eqref{poi1} to be a point $\hat x$ that satisfies:
\begin{equation}\label{asp} \hat x \in Q \ \ \ \mbox{and} \ \ \  \displaystyle\frac{f(\hat x) -f^*}{f^*-\ua} \ \le \ \varepsilon' \ . \end{equation}

Note that \eqref{asp} is a relative error measure, relative to the optimal bound gap $f^* - \ua$.  We focus on an $\varepsilon'$-relative solution rather than on an $\varepsilon$-absolute solution ($f(\hat x) \le f^* + \varepsilon$), as the former seems more natural in the setting where a strict lower bound is part of the problem description.  Indeed, consider the context of loss functions $f(\cdot)$ in statistics and machine learning where $\ua = 0$, in which case an $\varepsilon'$-relative solution $\hat x$ corresponds to $\frac{f(\hat x)}{f^*} \le (1+\varepsilon')$, and hence is a multiplicative measure of optimality tolerance.

Let $G$ denote the smallest scalar $\bar G$ satisfying:
\begin{equation}\label{gg} \dist (x,\optim )  \le \bar G \cdot (f(x) - \ua) \ \ \mbox{for~all~} x \in Q \ . \end{equation}
By its definition one sees that $G$ measures how fast the distances from the optimal solution set $\optim$ grow relative to the bound gap $f(x) - \ua$.  Therefore $G$ is a measure of the growth rate of the level sets of $f(\cdot)$.  We call $G$ the ``growth constant'' of the function $f(\cdot)$ for the given strict lower bound $\ua$.  Note that an equivalent definition of $G$ is given by:
\begin{equation}\label{ggg} G = \displaystyle\sup_{x \in Q} \left\{\frac{\dist (x,\optim )}{f(x) - \ua}\right\} \ . \end{equation}

Unlike the strict lower bound $\ua$, we do not assume that $G$ is known, nor do we need any upper bounds on $G$.  Indeed, neither knowledge of $G$ nor the finiteness of $G$ are needed in order to implement the computational methods presented herein; however the finiteness of $G$ is needed for the analysis of the methods to be meaningful.

We will see in Sections \ref{general} and \ref{mtc} that the knowledge of the fixed strict lower bound $\ua$ and the concept of the function growth constant $G$ lead to different versions of first-order methods with different computational guarantees than the traditional analysis of first-order methods would dictate.  Furthermore, these different computational guarantees can dominate the traditional guarantees in many cases but most notably when the initial iterate $x^0$ is far from the optimal solution. Roughly speaking, for several of the algorithms developed herein our computational guarantees grow like $\ln (1+\dist(x^0, \optim))$ in contrast to traditional guarantees where the growth is proportional to $\dist(x^0,\optim)$ and $\dist (x^0,\optim)^2$ (in the smooth and nonsmooth settings, respectively).


In a departure from typical optimization approaches to lower bounds such as those arising from duality theory wherein one desires as tight a lower bound as possible, herein the lower bound $\ua$ is strict, namely $\ua < f^*$, and it is fixed, i.e., it is not updated as part of a computational procedure.  It is best to think of this lower bound as a {\em structural} lower bound that is easily connected to known properties of the function $f(\cdot)$.    Such a strict lower bound on $f(\cdot)$ arises naturally in the settings of statistics and machine learning in the case of loss functions and/or regularization functions, see for example \cite{elements}.  Consider when $f(\cdot)$ is the logistic loss function $f(x) = \frac{1}{m}\sum_{i=1}^m \ln\left(1 + e^{-A_i x} \right)$ or the exponential loss function $f(x) = \ln\left(\sum_{i=1}^m e^{-A_i x} \right)$, perhaps with the addition of a regularization term $\lambda \| x\|_p^r$ for some $p\ge 1$, $r \ge 1$, and $\lambda \ge 0$.  If the sample data is not strictly separable, which translated herein means that there is no $x$ satisfying $Ax \ge 0$ unless $Ax=0$, then it follows that $f^* >0$ and so $\ua := 0$ is a strict lower bound and is quite natural in this setting.  Another example is regularized least-squares regression such as the LASSO and its cousins, wherein $f(\beta) = \tfrac{1}{2}\|\by - \bX\beta\|^2 + \lambda \| \beta\|_p^r$; it follows that $f^* \ge 0$ and one can assert that $f^* > 0 =: \ua$ under a variety of mild assumptions involving either $\lambda$ or the data matrix $\bX$.  Other classes of examples for which $\ua =0$ is a strict lower bound on $f^*$ include total variation (TV) loss functions which are used in image de-noising, as well as the broad class of minimum norm problems in general, under mild assumptions.  Another class of problems for which there is a natural strict lower bound on $f^*$ is the class of projectively transformed conic convex optimization problems under a particularly clever projective transformation, as developed by Renegar \cite{renegar2015}; indeed it was this problem class and the results in \cite{renegar2015} that gave rise to the line of research described herein.

We can interpret $G$ as connected to a lower estimator of $f(\cdot)$:  rearranging \eqref{gg}, we obtain:
\begin{equation}\label{gg2} f(x) \ \ge \ \bar f(x) := \ua + G^{-1}\dist(x,\optim) \ \ \mbox{for~all~} x \in Q \ . \end{equation}

\begin{figure}
\begin{center}
\includegraphics[scale=0.40]{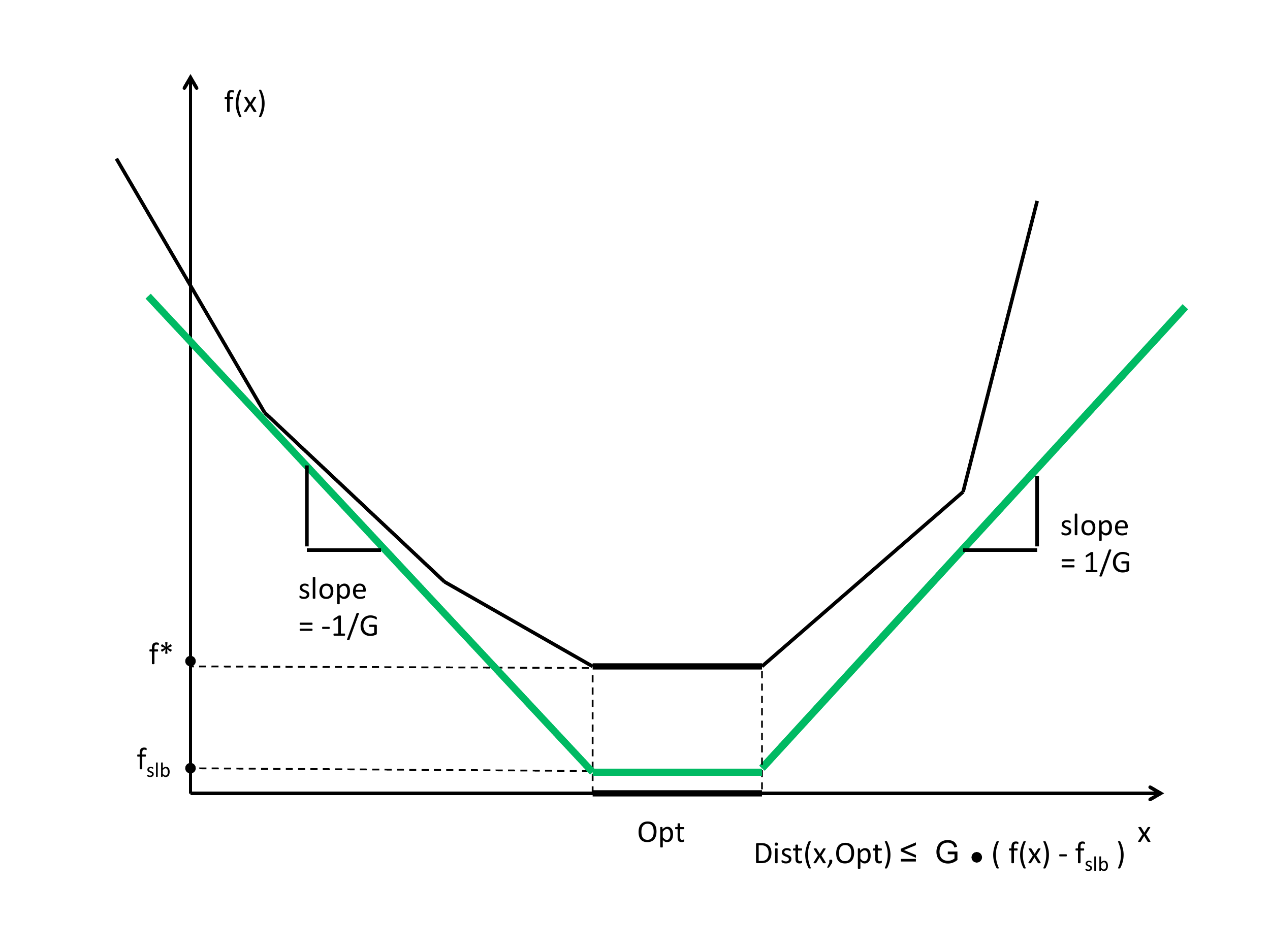}
\caption{Illustration of $G$ and $\ua$ for a function with multiple optimal solutions.}\label{fig1}
\end{center}
\end{figure}

Therefore the convex function $\bar f(x) = \ua + G^{-1}\dist(x,\optim)$ is a lower estimator of the function $f(\cdot)$ on $Q$.  This interpretation is illustrated in Figure \ref{fig1}.  As Figure \ref{fig1} illustrates, the concept of the growth constant $G$ is somewhat related to the notion of the modulus of weak sharp minima for \eqref{poi1}, see Polyak \cite{polyak1979sharp} and Burke and Ferris \cite{burke1993weak}; this relationship is discussed further in Appendix \ref{burke-ferris}.


A natural question to ask is under what circumstances is the growth constant $G$ finite?  Roughly speaking, it holds that $G$ is finite except when the objective function level sets are ill-behaved relative to their recession cone.  This is made precise in the following theorem, whose proof is given in Appendix \ref{app:lu-theorem}.  For $\varepsilon>0$, let $\optim_\varepsilon := \{ x \in Q :  \ f(x) \leq f^* + \varepsilon\}$ denote the $\varepsilon$-optimal level set of $f(\cdot)$ on $Q$, and let $S$ denote the recession cone of $\optim_\varepsilon$, namely $S:= \{ d \in \mathbb{R}^n : x+\theta d \in Q \ \mbox{and} \ f(x+\theta d) \le f^* + \varepsilon \ \mbox{for~all} \ \theta \ge 0 \}$.  Note that $S$ is the (common) recession cone of $\optim_\varepsilon$ for all $\varepsilon \ge 0$.\medskip

\begin{thm}\label{lu-theorem} Suppose that for some $\varepsilon>0$ there exists a bounded set $E_\varepsilon$ for which $\optim_\varepsilon \subset E_{\varepsilon} + S$ where $S$ is the recession cone of $\optim_\varepsilon$.  Then for any given strict lower bound $\ua < f^*$, the growth constant $G$ is finite. \qed
\end{thm}\medskip

Let us briefly examine special cases of Theorem \ref{lu-theorem}. Consider the case when $\optim = E + T$ where $E$ is a bounded convex set and $T$ is a subspace.  Then for any $\varepsilon > 0$ it is easy to show that $\optim_\varepsilon = E_{\varepsilon} + T$ for some bounded set $E_{\varepsilon}$, in which case Theorem \ref{lu-theorem} implies that $G$ is finite.  In particular, when $\optim$ itself is a bounded set, then we can set $T=\{0\}$, and so Theorem \ref{lu-theorem} implies that $G$ is finite.

For an example wherein $G=\infty$, consider the function $f(x_1,x_2):=\frac{x_2^2}{x_1}$ on $Q:=\{(x_1,x_2) : x_1 \ge 1\}$.  It is straightforward to check that the Hessian matrix $\nabla^2 f(x)$ is positive semidefinite on $Q$ and hence $f(\cdot)$ is convex on $Q$.  We have $f^*=0$ and $\optim = \left\{(x_1,0) : x_1 \ge 1 \right\}$.  However, the growth constant $G=\infty$ for any strict lower bound $\ua$, since by letting $(x_1,x_2) = (\beta^2, \beta)$ for any $\beta \ge 1$ we obtain using \eqref{ggg} that $$ G \ge \lim_{\beta \rightarrow + \infty} \frac{\dist \left((\beta^2, \beta),\optim\right)}{f\left(\beta^2,\beta\right)-\ua}=\lim_{\beta\rightarrow + \infty} \frac{\beta}{1-\ua} = +\infty \ . $$

\subsection{Overview of Results}

We use the knowledge of the fixed strict lower bound $\ua$ and the concept of the function growth constant $G$ to design and develop computational guarantees for new versions of first-order methods for solving the optimization problem \eqref{poi1}.  In Section \ref{general} we present such methods when $f(\cdot)$ is non-smooth and Lipschitz continuous with Lipschitz constant $M$.  In Theorem \ref{th1} we present an iteration complexity of $O\left(M^2G^2\left[\ln\left(1+\frac{f(x^0)-f^*}{f^*-\ua}\right) + \frac{1}{(\varepsilon')^2} \right]\right)$ for a version of Subgradient Descent that simultaneously runs with two step-sizes and occasional re-starting, which strictly improves the standard computational complexity bound for Subgradient Descent when $x^0$ is a ``cold start,'' i.e., $\dist(x^0, \optim)$ is large.  In the special case when the optimal objective function value $f^*$ is known,  Theorem \ref{th3} shows that the standard step-size rule for Subgradient Descent yields the same result.  And when $f(\cdot)$ can be smoothed, we present further improved computational guarantees for a new method (Algorithm \ref{smoothedscheme}) that successively smooths and restarts the Accelerated Gradient Method, see Theorem \ref{th4} herein.

In Section \ref{mtc} we present computational guarantees for new first-order methods when $f(\cdot)$ is smooth and has Lipschitz gradient with Lipschitz constant $L$.  We present a new first-order method (Algorithm \ref{smoothschemerestart}) based on periodically restarting the Accelerated Gradient Method, that leads to an iteration complexity of  $O\left(G\sqrt{L}\left[\sqrt{f(x^0)-\ua}  + \frac{\sqrt{f^*-\ua}}{\sqrt{\varepsilon'}} \right]\right)$ (Theorem \ref{th55}), which in many cases can improve the standard computational complexity bound for the Accelerated Gradient Method, most notably when $f(x^0)$ is far from the optimal value $f^*$ and $\varepsilon'$ is small.  And when $f(\cdot)$ has appropriate adjoint structure, we use parametric increased smoothing and restarting of the Accelerated Gradient Method to achieve a further improvement in the above computational guarantee (Theorem \ref{th5}).

Algorithm A in Renegar \cite{renegar2016a} provides an interesting approach to the general convex optimization setting, that bears comparison to the approach and results contained herein -- which are also designed for the general convex optimization setting.  Both Algorithm A in \cite{renegar2016a} and the algorithms herein generalize the methodology for conic optimization developed in Renegar \cite{renegar2014,renegar2015} to the general convex optimization problems, but they do so in different ways. Herein the generalization is obtained by introducing the new function measure $G$ based on the strict lower bound $\ua$, while in Algorithm A in \cite{renegar2016a} the original problem is transformed (implicitly or explicitly) to a conic optimization problem in a slightly lifted space. The resulting algorithms appear to be very different, and have different computational requirements and convergence bounds -- Algorithm A in \cite{renegar2016a} requires a $1$-dimensional root finding procedure each iteration, whereas Algorithm \ref{nsm} herein requires orthogonal projection onto the feasible region.  (And indeed it is rather remarkable that Algorithm A of \cite{renegar2016a} does not require such projection.)  Algorithm A does not need a Lipschitz constant; however in the case of a smooth objective function Algorithm A cannot take advantage of such smoothness, unlike Algorithm \ref{smoothschemerestart} (and also Algorithm \ref{smoothedscheme}) herein.  



The paper is organized as follows.  Section \ref{review} contains a brief review of the Subgradient Descent and an Accelerated Gradient Method.  Section \ref{general} contains first-order methods and computational guarantees when $f(\cdot)$ is non-smooth.  Section \ref{mtc} contains first-order methods and computational guarantees when $f(\cdot)$ is smooth.

\noindent {\bf Notation.}  Unless otherwise specified, the norm is the Euclidean (inner product) norm $\| x \| := \sqrt{x^Tx}$.  We occasionally refer to the $\ell_p$ norm of a vector $v$, which is denoted by $\|v\|_p$.
For $Q \subset \mathbb{R}^n$, let $\Pi_Q(\cdot)$ denote the Euclidean projection operator onto $Q$, namely $\Pi_Q(x) := \arg\min_{y \in Q} \|y-x\| $.  We define $\dist(x,S) := \min_y\{\|x-y\| : y \in S\}$.  The set of optimal solutions of \eqref{poi1} is denoted by $\optim := \{ x \in Q : f(x) =f^* \}$.
\section{Review of Subgradient Descent and an Accelerated Gradient Method}\label{review}
We briefly review the Subgradient Descent Method and an Accelerated Gradient Method (as analyzed in Tseng \cite{tseng}) for solving the convex optimization problem \eqref{poi1}.

\subsection{Subgradient Descent}  Recall that $g$ is a subgradient of $f(\cdot)$ at $x$ if the following subgradient inequality holds:  \begin{equation*}\label{jeff}f(y) \ge f(x) + g^T(y - x)  \ \ \ \mathrm{for~all~} y \in Q \ . \end{equation*} Let $\partial f(x)$ denote the set of subgradients of $f(\cdot)$ at $x$.  Here we assume that $f(\cdot)$ is Lipschitz continuous on a relatively open set $\hat{Q}$ containing $Q$, namely, there is a scalar $M$ for which
\begin{equation}\label{lip} |f(y) - f(x)| \le M \|y-x\| \ \ \ \ \mbox{for~all~} x,y\in \hat{Q} \ .
\end{equation}

It follows from \eqref{lip} that for all $x \in Q$ and $g \in \partial f(x)$ it holds that $\|g\| \le M$.

Algorithm \ref{subgradientmethod} presents the standard subgradient scheme.  In this method $x^k$ is the iterate at iteration $k$, the best objective value among the first $k$ iterates is $f^k_b$, and the best iterate among the first $k$ iterates is $x^k_b$.


\begin{algorithm}
\caption{Subgradient Method for Non-Smooth Optimization}\label{subgradientmethod}
$ \ $
\begin{algorithmic}
\STATE {\bf Initialize.} Initialize with $x^0 \in Q$, $f^0_b \gets f(x^0)$, $x^0_b \gets x^0$ . $ i \gets 0$ .  \\
$ \ $ \\
At iteration $i$ :\\
\STATE 1. {\bf Compute Subgradient.} Compute $g_i \in \partial f(x^i)$ .\\
\STATE 2. {\bf Determine Step-size.} Determine $\alpha_i \ge 0$ . \\
\STATE 3. {\bf Perform Updates.} Compute $x^{i+1} \gets \Pi_Q(x^i -\alpha_i g_i)$ , \\\medskip
\ \ \ \ \ \ \ \ \ \ \ \ \ \ \ \ \ \ \ \ \ \ \ \ \ \ \ \ \ \ \ \ $f^{i+1}_b\gets \min\{f^i_b, f(x^{i+1})\}$ , \\\medskip
\ \ \ \ \ \ \ \ \ \ \ \ \ \ \ \ \ \ \ \ \ \ \ \ \ \ \ \ \ \ \ \ $x^{i+1}_b \gets \argmin\limits_{x \in \{x^{i}_b, \ x^{i+1} \}} \{f(x)\}$ . \\
$ \ $ \\
\end{algorithmic}
\end{algorithm}

The following theorem summarize well-known computational guarantees associated with the subgradient descent method.\medskip

\begin{thm}{\bf{(Convergence Bounds for Subgradient Descent}\ \cite{polyakbook, nesterovBook})}\label{subgrad}\\
\noindent (i) Consider the subgradient descent method (Algorithm \ref{subgradientmethod}). Then for all $k \ge i \geq 0$, the following inequality holds:
\begin{equation*}\label{sd_bound1}
f^k_b \ \ \leq \ \ f^* \ + \ \frac{\dist(x^i, \optim)^2+ \sum_{l=i}^k\|g_l\|^2\alpha_l^2}{2\sum_{l=i}^k\alpha_l}  \ \ \leq \ \ f^* \ + \ \frac{\dist(x^i, \optim)^2+ M^2\sum_{l=i}^k\alpha_l^2}{2\sum_{l=i}^k\alpha_l} \ .
\end{equation*}

\noindent (ii) Suppose that $f^*$ is known, and let the step-sizes for Algorithm \ref{subgradientmethod} be $\alpha_i = (f(x^i) - f^*)/\|g_i\|^2$.  Then for all $k \ge i \geq 0$, the following inequality holds:
\begin{equation*}\label{sd_bound2}
\ \ \ \ \ \ \ \ \ \ \ \ \ \ f^k_b    \ \ \leq \ \ f^* \ + \ \frac{M \dist(x^i, \optim)}{\sqrt{k-i+1}} \ . \ \ \ \ \ \ \ \ \ \ \ \ \ \ \ \ \ \ \ \ \ \qed
\end{equation*}
\end{thm}

Suppose that we seek to bound the number of iterations $N$ of the Subgradient Descent method required to compute an (absolute) $\varepsilon$-optimal solution of \eqref{poi1}, which is a point $\hat x \in Q$ that satisfies $f(\hat x) \le f^* + \varepsilon$.  If $\varepsilon > 0$ is given, and the step-sizes are chosen as $\alpha_i = \varepsilon/\|g_i\|^2$, then it follows from part $(i)$ of Theorem \ref{subgrad} that $f^N_b \le f^* + \varepsilon$ for all
\begin{equation}\label{friday1}
N \ge \bar N := \displaystyle\frac{M^2\dist(x^0, \optim)^2}{\varepsilon^2} -1 \ .
\end{equation}  If instead we know (or can bound from above) $\dist(x^0, \optim)$, and the step-sizes are chosen as $\alpha_i = \dist(x^0, \optim)/(\sqrt{N+1}\|g_i\|)$ where $N$ satisfies \eqref{friday1}, then it also follows from part $(ii)$ of Theorem \ref{subgrad} that $f^N_b \le f^* + \varepsilon$.  And if $f^*$ is known, then the bound \eqref{friday1} is also sufficient to guarantee $f^N_b \le f^* + \varepsilon$ if the steps-sizes are chosen as in part $(ii)$ of Theorem \ref{subgrad}.  Furthermore, it follows from \cite{NemirovskyYudin83} that the bound \eqref{friday1} cannot in general be improved in the black-box oracle model of computation with complexity bounds depending only on $M$, $\dist(x^0,\optim)$, and $\varepsilon$.  In this regard, we note that the dependence on additional parameters, namely the strict lower bound $\ua$ and the function growth constant $G$, which are used throughout this paper, shows how we can achieve different (and better in many cases) complexity bounds by including additional parameters and appropriately amending algorithms and their analysis.

\subsection{Accelerated Gradient Method for Smooth Optimization}  Here we assume that $f(\cdot)$ is differentiable with Lipschitz continuous gradient on $Q$, namely, there is a scalar $L$ for which
\begin{equation}\label{lipgrad} \|\nabla f(y) - \nabla f(x)\| \le L \|y-x\| \ \ \ \ \mbox{for~all~} x,y\in Q \ .
\end{equation}
Algorithm \ref{accsimplegradientscheme} presents a standard Accelerated Gradient Method as in Tseng \cite{tseng}.   

\begin{algorithm}
\caption{Accelerated Gradient Method}\label{accsimplegradientscheme}
$ \ $
\begin{algorithmic}
\STATE {\bf Initialize.}  Initialize with $x^0 \in Q$ and $z^0 := x^0$, and $i \gets 0$ .  Define step-size parameters $\theta_i \in (0,1]$ recursively by $\theta_0 :=1$ and $\theta_{i+1}$ satisfies $\tfrac{1}{\theta_{i+1}^2} - \tfrac{1}{\theta_{i+1}} = \tfrac{1}{\theta_i^2}$ . \\
$ \ $ \\
At iteration $i$ :\\
\STATE 1. {\bf Perform Updates.}  Define $y^i \gets (1-\theta_i)x^i + \theta_i z^i$, and compute $\nabla f(y^i)$ , \\\medskip
\ \ \ \ \ \ \ \ \ \ \ \ \ \ \ \ \ \ \ \ \ \ \ \ \ \ \ \ \ \ \ \ $z^{i+1} \gets \arg\min_{x \in Q} \{f(y^i) + \nabla f(y^i)^T(x-z^i) + \tfrac{1}{2}\theta_iL(x-z^i)^T(x-z^i)\}$ ,\\\medskip
\ \ \ \ \ \ \ \ \ \ \ \ \ \ \ \ \ \ \ \ \ \ \ \ \ \ \ \ \ \ \ \ $x^{i+1} \gets (1-\theta_i)x^i + \theta_i z^{i+1}$ .\\\medskip
\end{algorithmic}
\end{algorithm}

For $\delta \ge f^*$ define the level set $S_\delta := \{ x \in Q : f(x) \le  \delta \}$.  For $x \in Q$, let $ \dist (x, S_\delta )$ denote the distance from $x$ to the level set $S_\delta$, namely $ \dist (x,S_
\delta ) := \min_y\{\|y-x\| : y \in S_\delta\}$.  The following theorem is a computational guarantee for the Accelerated Gradient Method due to Tseng \cite{tseng}.\medskip

\begin{thm}{\bf{(Convergence Bound for Accelerated Gradient Method}\ \cite{tseng})}\label{accgrad}
Consider the Accelerated Gradient Method (Algorithm \ref{accsimplegradientscheme}). Let $\delta \ge f^*$ and $S_\delta :=\{x\in Q : f(x) \le \delta \}$.  Then for all $k \ge 0$, the following inequality holds:
\begin{equation*}\label{sd_bound3}
\ \ \ \ \ \ \ \ \ \ \ \ \ \ f(x^k)  -\delta  \ \ \leq \ \  \frac{2L \dist(x^0, S_\delta)^2}{(k+1)^2} \ . \ \ \ \ \ \ \ \ \ \ \ \ \ \ \qed
\end{equation*}
\end{thm}

Note that in the case when $\delta = f^*$, then $S_\delta = \optim$ whereby Theorem \ref{accgrad} specializes to the standard result for the Accelerated Gradient Method.  We will utilize the more general result in Theorem \ref{accgrad} in the context of smoothing of a non-smooth function, in Sections \ref{generalissimo} and \ref{mtc} herein.

\section{Computational Guarantees when $f(\cdot)$ is Non-Smooth}\label{general}

Let $\varepsilon' > 0$ be given.  We aspire to compute an $\varepsilon'$-relative solution of \eqref{poi1}, which recall from \eqref{asp} is a point $\hat x \in Q$ satisfying:
$\frac{f(\hat x) -f^*}{f^*-\ua} \ \le \ \varepsilon'$.  In this section we present three new computational guarantees for first-order methods applied to computing a $\varepsilon'$-relative solution of problem \eqref{poi1} that are based on the strict lower bound $\ua$ and growth constant $G$.  The first guarantee is for a new algorithm based on Subgradient Descent that runs two different step-sizes simultaneously with occasional re-starts.  The second guarantee is for the standard Subgradient Method using a standard step-size rule in the case when the optimal value $f^*$ is known.  The third guarantee is for the case when the function $f(\cdot)$ can be smoothed and then solved using an algorithm based on the Accelerated Gradient Method.

\subsection{Subgradient Descent using Two Step-Size Rules Running Simultaneously}\label{sbs}

We consider solving \eqref{poi1} using a version of subgradient descent that simultaneously runs two versions of the Subgradient Descent Method -- each with a different step-size rule -- with occasional simultaneous re-starts of both versions.  The formal description of our method is given in Algorithm \ref{nsm}.  In the algorithm, the notation ``$(x_{i,j+1}, f^{i,j+1}_b, x^{i,j+1}_b) \gets \mathrm{SDM}(f(\cdot),x_{i,j}, \alpha_{ij}, g_{ij})$'' denotes assigning to $x_{i,j+1}$ the next value of the Subgradient Descent Method applied to the optimization problem \eqref{poi1} with objective function $f(\cdot)$ with current point $x_{i,j} \in Q$ using the step-size $\alpha_{ij}$ and the subgradient $g_{ij}$, along with updates of the best objective function value obtained thus far $f^{i,j+1}_b$ with the corresponding best iterate computed $x^{i,j+1}_b$.

\begin{algorithm}[hbt]
\caption{Non-Smooth Method with Two Step-Size Rules Running Simultaneously}\label{nsm}
$ \ $
\begin{algorithmic}
\STATE {\bf Initialize.}  Initialize with $x^0 \in Q$ and $\varepsilon' >0$ .  \\

Define constants $\bar\varepsilon' := 0.9$ ,  $\varepsilon := \tfrac{\varepsilon'}{1+\varepsilon'}$ , $\bar\varepsilon := \tfrac{\bar\varepsilon'}{1+\bar\varepsilon'}$, $B := 1/\sqrt{e}$ , $F := \sqrt{e}$ . \\

Set $x_{1,0} \gets x^0$ , $\bar x_{1,0} \gets x^0$, $ i \gets 1$  .
$ \ $ \\
$ \ $ \\
At outer iteration $i$ :

\STATE 1. {\bf Initialize inner iterations.} $f^{i,0}_b \gets f(x_{i,0})$ \ , \ $\bar f^{i,0}_b \gets f(\bar x_{i,0})$   \\\medskip
\ \ \ \ \ \ \ \ \ \ \ \ \ \ \ \ \ \ \ \ \ \ \ \ \ \ \ \ \ \ \ \ \ \ \ \ \ \ \ \ \ \ \ \ $x^{i,0}_b \gets x_{i,0}$ \ \ \ \ \ , \ $\bar x^{i,0}_b \gets \bar x_{i,0}$ \\\medskip
\ \ \ \ \ \ \ \ \ \ \ \ \ \ \ \ \ \ \ \ \ \ \ \ \ \ \ \ \ \ \ \ \ \ \ \ \ \ \ \ \ \ \ \ $K_i \gets +\infty$ , $j \gets 0$ . \\\medskip
\STATE 2. {\bf Test/update current iterates.} At inner iteration $j$: \\
 \ \ \ \ \ \ \ \ \ \ \ \ \ If \ \ $\displaystyle\frac{f(x_{i,j})-\ua }{f(x_{i,0})-\ua} \ge B$ \ \ and \ \ $\displaystyle\frac{f(\bar x_{i,j})-\ua }{f(x_{i,0})-\ua} \ge B$, then \\\medskip
\ \ \ \ \ \ \ \ \ \ \ \ \ \ \ \ \ \ \ \ \ \ \ \ \ {\bf Compute subgradients.} Compute $g_{ij} \in \partial f(x_{i,j})$ , $\bar g_{ij} \in \partial f(\bar x_{i,j})$  \\\medskip
\ \ \ \ \ \ \ \ \ \ \ \ \ \ \ \ \ \ \ \ \ \ \ \ \ {\bf Set step-sizes.} $\alpha_{ij} \gets \tfrac{\varepsilon (f(x_{i,0})-\ua)}{F\|g_{ij}\|^2}$ \ \ \ \ \ \ \ \ \  , \ \ $\bar\alpha_{ij} \gets \tfrac{\bar\varepsilon (f(\bar x_{i,0})-\ua)}{F\|\bar g_{ij}\|^2}$ \\\medskip
\ \ \ \ \ \ \ \ \ \ \ \ \ \ \ \ \ \ \ \ \ \ \ \ \ {\bf Update:}  \ \ \ \ \ \ $(x_{i,j+1}, f^{i,j+1}_b, x^{i,j+1}_b) \gets \mathrm{SDM}(f(\cdot),x_{i,j}, \alpha_{ij}, g_{ij})$\\\medskip
\ \ \ \ \ \ \ \ \ \ \ \ \ \ \ \ \ \ \ \ \ \ \ \ \ \ \ \ \ \ \ \ \  \ \ \ \ \ \ \ \ \ \ \ $(\bar x_{i,j+1}, \bar f^{i,j+1}_b, \bar x^{i,j+1}_b) \gets \mathrm{SDM}(f(\cdot), \bar x_{i,j}, \bar \alpha_{ij}, \bar g_{ij})$\\\medskip
 \ \ \ \ \ \ \ \ \ \ \ \ \ Else if \ \ $\displaystyle\frac{f(x_{i,j})-\ua }{f(x_{i,0})-\ua} < B$ , then: \\\medskip
  \ \ \ \ \ \ \ \ \ \ \ \ \ \ \ \ \ \ \ \ \ \ \ \  \ \ \ \ \ \ \ \ \ \ \ \ \ \ \ \ $K_i \gets j$ , $x_{i+1,0} \gets x_{i,j}$ , $\bar x_{i+1,0} \gets x_{i,j}$ , $i \gets i+1$ ,  and Goto Step 1.\\\medskip
 \ \ \ \ \ \ \ \ \ \ \ \ \ Else \ \ \ \ \ $\displaystyle\frac{f(\bar x_{i,j})-\ua }{f(\bar x_{i,0})-\ua} < B$ , and: \\\medskip
  \ \ \ \ \ \ \ \ \ \ \ \ \ \ \ \ \ \ \ \ \ \ \ \  \ \ \ \ \ \ \ \ \ \ \ \ \ \ \ \ $K_i \gets j$ , $x_{i+1,0} \gets \bar x_{i,j}$ , $\bar x_{i+1,0} \gets \bar x_{i,j}$ , $i \gets i+1$ ,  and Goto Step 1.\\\medskip

\end{algorithmic}
\end{algorithm}\medskip

We now walk through the structure of Algorithm \ref{nsm}.  The algorithm requires as input the starting point $x^0$ and the desired relative accuracy value $\varepsilon'$ used to define an $\varepsilon'$-relative solution, see \eqref{asp}.   The algorithm then defines an absolute constant $\bar\varepsilon' := 0.9$.  The two values $\varepsilon'$ and $\bar\varepsilon'$ are then used as aspirational goals for simultaneously running the standard Subgradient Descent Method in search of either an $\varepsilon'$-relative solution of \eqref{poi1} or an $\bar \varepsilon'$-relative solution of \eqref{poi1}.  For notational ease, both $\varepsilon'$ and $\bar\varepsilon'$ are converted to a slightly different form by defining $\varepsilon$ and $\bar\varepsilon$.  At the start of the $i^{\mathrm{th}}$ outer iteration, Algorithm \ref{nsm} runs the Subgradient Descent Method simultaneously using two different step-size rules (but starting at the same point $x_{i,0}=\bar x_{i,0}$), and so generates inner iterations $\{x_{i,j}\}$ and $\{\bar x_{i,j}\}$ for $j=0,1, \ldots$ based on computed subgradients $\{g_{ij}\}$ and $\{\bar g_{ij}\}$ and step-sizes $\{\alpha_{ij}\}$ and $\{\bar\alpha_{ij}\}$, respectively.  The only structural difference between the two instantiations of Subgradient Descent is that the steps-sizes $\{\alpha_{ij}\}$ use $\varepsilon$ in their definition whereas $\{\bar\alpha_{ij}\}$ use $\bar\varepsilon$ in their definition.  The number of inner iterations $j$ that are run in the  $i^{\mathrm{th}}$ outer iteration is initially set to be $K_i \gets +\infty$.  If either $f(x_{i,j})$ or $f(\bar x_{i,j})$ makes sufficient progress relative to the starting value $f(x_{i,0}) (=f(\bar x_{i,0}))$ as determined in the ratio test at the start of Step (2.), then the outer iteration $i$ is concluded and $K_i$, which counts the number of inner iterations therein, is updated to $K_i \gets j$.  Finally, the next outer iteration starting values $x_{i+1,0} = \bar x_{i+1,0}$ are re-set to either  $x_{i,j}$ or $\bar x_{i,j}$ , depending on which of $x_{i,j}$ or $\bar x_{i,j}$ satisfies the ratio test.

Many of the ideas used in the construction of Algorithm \ref{nsm} were motivated from similar notions developed in Algorithm 2 of \cite{renegar2015} as well as the algorithm ``MainAlgo'' in \cite{renegar2016b} (which uses the construct of running two algorithms simultaneously with different parameters).

Regarding counting of iterates $x_{i,j}$, $\bar x_{i,j}$ that are computed of Algorithm \ref{nsm}, we will say that the algorithm has computed an iterate whenever it computes a subgradient and then calls $\mathrm{SDM}( \cdot, \cdot, \cdot)$.  There are therefore two iterates computed at each inner iteration.  We have: \medskip

\begin{thm}{\bf (Complexity Bound for Algorithm \ref{nsm})}\label{th1} Within a total number of iterates computed that does not exceed $$ 18M^2 G^2 \left( 2.7 \ln\left( 1 + \frac{f(x^0)-f^*}{f^* - \ua} \right) \ +  \  \left(\frac{1+\varepsilon'}{\varepsilon'}\right)^2 \right) \ , $$
Algorithm \ref{nsm} will compute an iterate $x_{i,j}$ for which $$ \frac{f(x_{i,j}) - f^*}{f^* - \ua} \le \varepsilon'  \  . \ \qed $$
\end{thm} \medskip

Since $f(x^0) \le f^* + M\dist(x^0,\optim)$, the computational guarantee in Theorem \ref{th1} can itself be bounded by:
\begin{equation}\label{jim}18M^2 G^2 \left( 2.7 \ln\left( 1+ \frac{M\dist(x^0,\optim)}{f^*- \ua}\right) \ +  \  \left(\frac{1+\varepsilon'}{\varepsilon'}\right)^2 \right) \ , \end{equation}

which is qualitatively different from the guarantee of the standard Subgradient Descent Method (Algorithm \ref{subgradientmethod}) in \eqref{friday1} in two interesting ways.  First, the dependence in \eqref{friday1} on $\dist(x^0,\optim)$ is quadratic, whereas in \eqref{jim} it is logarithmic.  Second, although both guarantees are linear in the inverse square of the desired relative accuracy $\varepsilon'$ (from \eqref{asp} an $\varepsilon'$-relative solution corresponds to an absolute $\varepsilon' \cdot \left(f^*-\ua\right)$ solution of \eqref{poi1}), however $x^0$ affects this factor multiplicatively through $\dist(x^0,\optim)^2$ in \eqref{friday1}, whereas the factor is independent of $x^0$ in \eqref{jim}.

Let us also quantitatively compare the computational guarantee of Theorem \ref{th1} with the standard guarantee for Subgradient Descent given by \eqref{friday1}.  The standard computational guarantee \eqref{friday1} can be written as:
\begin{equation*}
\frac{M^2\dist (x^0,\optim)^2}{\varepsilon'^2\left(f^*-\ua\right)^2} \ .
\end{equation*}
Let us presume that $\varepsilon'$ is small, whereby $\frac{1+\varepsilon'}{\varepsilon'}\approx \frac{1}{\varepsilon'}$.  Then the ratio of the new guarantee \eqref{jim} from Theorem \ref{th1} to the standard guarantee \eqref{friday1} is at most \begin{equation}\label{wow} \frac{\mathrm{Guarantee~of~Theorem~\ref{th1}}}{\mathrm{Standard~Guarantee~\eqref{friday1}}} \ \le \ 18(f^*-\ua)^2G^2\left(\frac{2.7 (\varepsilon')^2 \ln \left(1+ \frac{M\dist(x^0,\optim)}{f^*- \ua}\right) + 1 }{\dist(x^0, \optim)^2} \right) \ .  \end{equation}
Notice from \eqref{wow} that for any instance of \eqref{poi1}, when $\dist(x^0, \optim)$ is sufficiently large the right-hand side of \eqref{wow} can be made arbitrarily small, thereby showing that in these cases the computational guarantee in Theorem \ref{th1} can be made arbitrarily better than the standard guarantee \eqref{friday1} for  Subgradient Descent.

We will prove Theorem \ref{th1} by first establishing eight propositions.  The reader familiar with \cite{renegar2015} will notice certain resemblances between aspects of the proof constructs below and the proof of Theorem 3.8 of \cite{renegar2015}, see also \cite{renegar2016a}.  Throughout, for notational convenience, we will work with three constants $B $, $F$, and $\bar\varepsilon' $ that must be chosen to satisfy the conditions:
\begin{equation*}\label{conditions}
B \in (0,1) \ \ , \ \ F > \frac{1}{2B} \ \ , \ \ \mbox{and} \ \ \bar \varepsilon' > 0  \ ,
\end{equation*} and whose specific values in Algorithm \ref{nsm} are set to $B = 1/\sqrt{e}$, $F = \sqrt{e}$, and $\bar\varepsilon' = 0.9$, where $e$ is the base of the natural logarithm.\medskip

Let $\delta'>0$ play the role of either $\varepsilon'$ or $\bar \varepsilon'$, and also define $\delta := \frac{\delta'}{1+\delta'}$ (analogous to the definitions of $\varepsilon$ and $\bar \varepsilon$).\medskip

The first two propositions below apply to the generic setting of the Subgradient Descent Method. \medskip

\begin{prop}\label{prop1}  Let $\delta \in (0,1)$ be given, and suppose we run the Subgradient Descent Method (Algorithm \ref{subgradientmethod}) with starting iterate $\hat x^0$, using step-sizes: $$ \alpha_{j} := \frac{\delta (f(\hat x^0) - \ua)}{F\|g_{j}\|^2} $$ for all iterations $j$.  Then for all $j \ge 0$ it holds that
$$f^{j}_b - \ua \le f^* - \ua + \left[ \frac{G^2M^2F}{2\delta(j+1)}+ \frac{\delta}{2F} \right](f(\hat x^0) - \ua) \ . $$
\end{prop}

{\bf Proof:}  Define $\alpha:= \frac{\delta(f(\hat x^0) - \ua)}{F}$.  Then $\alpha_j = \frac{\alpha}{\|g_j\|^2} \ge \frac{\alpha}{M^2}$.  It follows from part $(i)$ of Theorem \ref{subgrad} that
\begin{equation*}\label{101}\begin{array}{rcl}
f^{j}_b - \ua & \le & f^* - \ua + \displaystyle\frac{\dist(\hat x^0, \optim)^2}{2\sum_{l=0}^{j}\alpha_l}+ \displaystyle\frac{ \sum_{l=0}^{j}\|g_l\|^2\alpha_l^2}{2\sum_{l=0}^{j}\alpha_l} \\ \\
& \le & f^* - \ua + \displaystyle\frac{M^2\dist(\hat x^0, \optim)^2}{2\alpha(j+1)}+ \displaystyle\frac{ \alpha}{2} \\ \\
& \le & f^* - \ua + \displaystyle\frac{M^2G^2F(f(\hat x^0) - \ua)^2}{2(j+1)\delta(f(\hat x^0) - \ua)}+ \displaystyle\frac{ \delta(f(\hat x^0) - \ua)}{2F} \   , \end{array}
\end{equation*}
where the second inequality uses the definition of $\alpha$ and the inequality $\|g_j\| \le M$, and the third inequality uses the definition of $G$.  Simplifying the last expression completes the proof. \qed \medskip

\begin{prop}\label{prop2}  Under the identical set-up as Proposition \ref{prop1}, let $\delta' := \delta/(1-\delta)$, and define:
\begin{equation*}
W := \left\lfloor \frac{FM^2G^2}{2\delta^2\left[ B-\tfrac{1}{2F} \right]} \right\rfloor \ .
\end{equation*}  Then either $\frac{f^{W}_b -f^*}{f^*-\ua} \le \delta'$ ,  or $f^{W}_b - \ua \le B(f(\hat x^0) -\ua) $, or both.
\end{prop}
{\bf Proof:}  Suppose that $\frac{f^{W}_b -f^*}{f^*-\ua} > \delta'$.  This rearranges to:
$\delta' <  \frac{f^{W}_b - \ua }{f^*-\ua} -1$, whereby 
\begin{equation}\label{907}
 \frac{f^*-\ua}{f^{W}_b - \ua} < \frac{1}{1+\delta' } = 1-\delta \ .
\end{equation}
Invoking Proposition \ref{prop1} we have:
\begin{equation*}\begin{array}{rcl} f^{W}_b - \ua  &\le&   f^* - \ua + \left[ \frac{G^2M^2F}{2\delta(W+1)}+ \frac{\delta}{2F} \right](f(\hat x^0) - \ua) \\ \\
& < &   f^* - \ua + \left[ \delta \left(B - \frac{1}{2F}\right) + \frac{\delta}{2F}\right](f(\hat x^0) - \ua) \\ \\
& = & f^* - \ua + \delta B (f(\hat x^0) - \ua) \\ \\
& < & (1-\delta)(f^{W}_b - \ua) + \delta B (f(\hat x^0) - \ua) \  ,
\end{array}\end{equation*}
where the second inequality follows since $W+1 > \frac{FM^2G^2}{2\delta^2[B-1/(2F)]}$, and the last inequality uses \eqref{907}.  Rearranging the final inequality and dividing by $\delta$ then yields $f^{W}_b -\ua  \le B(f(\hat x^0) -\ua) $, which completes the proof. \qed \medskip

In the next two propositions we apply Proposition \ref{prop2} directly to the setting of Algorithm \ref{nsm}. \medskip

\begin{prop}\label{prop3}  Consider outer iteration $i$ of Algorithm \ref{nsm}. Define:
\begin{equation*}
U := \left\lfloor \frac{FM^2G^2}{2\varepsilon^2\left[ B-\tfrac{1}{2F} \right]} \right\rfloor \ .
\end{equation*}  If $K_i > U$, then $\frac{f^{ij}_b -f^*}{f^*-\ua} \le \varepsilon'$ for all $j = U, \ldots, K_i$ .
\end{prop}

{\bf Proof:}  Let us apply Proposition \ref{prop2} with $\delta':=\varepsilon'$, $W:=U$, and $\hat x^0 := x_{i,0}$.  If $K_i >U$, then by definition of $K_i$ it holds that $f^{i,U}_b - \ua > B(f^{i,0}_b - \ua)$.  Therefore from Proposition \ref{prop2} it holds that
$\frac{f^{i,U}_b -f^*}{f^*-\ua} \le \varepsilon'$.  \qed \medskip

\begin{prop}\label{prop4}  Consider outer iteration $i$ of Algorithm \ref{nsm}. Define:
\begin{equation*}
V := \left\lfloor \frac{FM^2G^2}{2\bar\varepsilon^2\left[ B-\tfrac{1}{2F} \right]} \right\rfloor \ .
\end{equation*}  If $K_i > V$, then $\frac{f(x_{i,0}) -\ua}{f^*-\ua} \le \frac{1+\bar\varepsilon'}{B}$ .
\end{prop}

{\bf Proof:}  Let us similarly apply Proposition \ref{prop2} with $\delta':=\bar\varepsilon'$, $W:=V$, and $\hat x^0 := x_{i,0}$.  If $K_i >V$, then by definition of $K_i$ it holds that $f^{i,V}_b - \ua > B(f(x_{i,0}) - \ua)$.  Therefore from Proposition \ref{prop2} it holds that
$\frac{f^{i,V}_b -f^*}{f^*-\ua} \le \bar\varepsilon'$.  Combining these inequalities we obtain:
\begin{equation*}
B(f(x_{i,0}) - \ua) < f^{i,V}_b - \ua = f^{i,V}_b - f^* + f^* - \ua \le \bar\varepsilon'(f^*-\ua) + f^* - \ua \ ,
\end{equation*} and rearranging yields the result.  \qed \medskip

In the next proposition we use the standard notation $ a^+$ for the nonnegative part of a scalar $a$.

\begin{prop}\label{prop5}  Let $m$ denote the number of outer iterations $i$ of Algorithm \ref{nsm} for which
\begin{equation*}
\frac{f(x_{i,0}) - \ua}{f^* - \ua} \ > \ \frac{1+\bar\varepsilon'}{B} \ .
\end{equation*}Then
\begin{equation*}\label{ce}
m \ \le \  {\left\lceil \frac{\ln\left(1+\frac{f(x^0)-f^*}{f^* - \ua}\right)- \ln \left( \frac{1+\bar\varepsilon'}{B} \right)}{\ln(1/B)} \right\rceil}^+ \ .
\end{equation*}
\end{prop}

{\bf Proof:}  If $m=0$ then the result holds trivially, so let us suppose that $m \ge 1$.  It then follows using induction on $f(x_{i+1, 0}) - \ua \le B(f(x_{i, 0}) - \ua)$ that
$$\frac{1+\bar\varepsilon'}{B}< \frac{f(x_{m,0}) - \ua}{f^* - \ua} \le \frac{B^{m-1}(f(x_{1,0}) - \ua)}{f^* - \ua}  \ , $$ and taking logarithms yields $$m-1 \ < \ \frac{\ln\left(\displaystyle\tfrac{f(x_{1,0})-\ua}{f^* - \ua}\right) - \ln \left( \frac{1+\bar\varepsilon'}{B} \right)}{\ln(1/B)} \ = \ \frac{\ln\left(1 + \displaystyle\tfrac{f(x^0)-f^*}{f^* - \ua}\right) - \ln \left( \frac{1+\bar\varepsilon'}{B} \right)}{\ln(1/B)}   \ , $$ from which the result follows. \qed \medskip

In the following proposition, as well as others later on, we use the standard notational convention that $\sum_{i=1}^n \cdot := 0 $ for $n \le 0$. \medskip

\begin{prop}\label{prop5a} Let $V$ and $m$ be as defined in Propositions \ref{prop4} and \ref{prop5}.  Then $x_{m+1,0}$ exists, and and let $T_m$ denote the total number of iterates computed prior to and including $x_{m+1,0}$.   It holds that: $$ \frac{f(x_{m+1,0}) - \ua}{f^* - \ua} \le \frac{1+\bar\varepsilon'}{B}  \ , $$
and furthermore $T_m \le 2mV$ .  \end{prop}

{\bf Proof:} If $m=0$ then the results holds trivially since
$$ \frac{f(x_{1,0}) - \ua}{f^* - \ua} \le \frac{1+\bar\varepsilon'}{B}  \ . $$
Next suppose that $m \ge 1$, and consider any outer iteration $i \le m$.  Then since
\begin{equation*}
\frac{f(x_{i,0}) - \ua}{f^* - \ua} \ > \ \frac{1+\bar\varepsilon'}{B} \ ,
\end{equation*} it follows from Proposition \ref{prop4} that $K_i \le V$.  This also implies that $x_{m+1,0}$ exists and therefore must satisfy
$$ \frac{f(x_{m+1,0}) - \ua}{f^* - \ua} \le \frac{1+\bar\varepsilon'}{B}  \ . $$  Finally, since $T_m = 2\sum_{i=1}^m K_i$, it therefore follows that $T_m \le 2mV$ . \qed \medskip

\begin{prop}\label{prop6}  Let $p$ denote the number of outer iterations $i$ for which $K_i$ is finite.  Then
\begin{equation*}
p \le m + \left\lfloor\frac{\ln \left( \frac{1+\bar\varepsilon'}{B} \right)}{\ln(1/B)}\right\rfloor \ ,
\end{equation*} where $m$ is as defined in Proposition \ref{prop5}.
\end{prop}

{\bf Proof:}  It follows from Proposition \ref{prop5a} that $p \ge m$.  Therefore $f^* - \ua \le f(x_{p,K_p}) - \ua =  f(x_{p+1,0}) - \ua \le B^{p-m}(f(x_{m+1,0})-\ua) \le B^{p-m}\left(\frac{1+\bar\varepsilon'}{B}\right) (f^*-\ua)$, where we have used the properties of $x_{m+1,0}$ in Proposition \ref{prop5a}.  Taking logarithms yields $$p-m \  \le \ \frac{\ln\left(\frac{1+\bar\varepsilon'}{B}\right)}{\ln(1/B)}  \ , $$ from which the result follows. \qed \medskip

\begin{prop}\label{prop6a}  Let $U$, $m$, and $p$ be as defined in Propositions \ref{prop3}, \ref{prop5}, and \ref{prop6}.  Within a total number of computed iterates after $x_{m+1,0}$ that does not exceed
$2(p-m+1)U$,  Algorithm \ref{nsm} will compute an iterate $x_{\hat i, \hat j}$ for which $$ \frac{f(x_{\hat i, \hat j}) - f^*}{f^* - \ua} \le \varepsilon'  \  . $$   \end{prop}

{\bf Proof:}  Let $\hat i$ denote the index of the first outer iteration $i \in \{m+1, \ldots, p+1\}$ for which $K_i > U$.  Notice that since $K_{p+1} = +\infty$ it must hold that $\hat i \le p+1$.  It follows from Proposition \ref{prop3} that $\frac{f^{\hat i,U}_b -f^*}{f^*-\ua} \le \varepsilon'$ and hence for some $\hat j \le U$ it holds that $\frac{f(x_{\hat i,\hat j}) -f^*}{f^*-\ua} \le \varepsilon'$.  Let us now count the number of iterates computed after $x_{m+1,0}$ and prior to and including $x_{\hat i, \hat j}$.  This number is bounded above by:
$$2\left(\sum_{i=m+1}^{\hat i -1}K_i + U \right) \le  2\left((\hat i -m-1)U + U \right) = 2(\hat i -m)U \le 2(p-m+1)U  \ , $$
where the first inequality follows since $K_i \le U$ for $i < \hat i$, and the last inequality uses $\hat i \le p+1$. \qed \medskip

We now use these propositions to prove Theorem \ref{th1}.\medskip

{\bf Proof of Theorem \ref{th1}:} Utilizing the definitions of $U$, $V$, $m$, $p$, and $x_{\hat i, \hat j}$ in Propositions \ref{prop3}, \ref{prop4}, \ref{prop5}, \ref{prop6}, and \ref{prop6a}, it follows from Propositions \ref{prop5a} and \ref{prop6a} that the total number of iterates computed prior to and including $x_{\hat i, \hat j}$ is at most $2[mV + (p-m+1)U]$.  Substituting the values of $U$ and $V$ and using the bounds on $m$ and $p$ in Propositions \ref{prop5} and \ref{prop6} yields:
$$\begin{array}{ll}
 & 2{\left\lceil \frac{\ln\left(1+\frac{f(x^0)-f^*}{f^* - \ua}\right)- \ln \left( \frac{1+\bar\varepsilon'}{B} \right)}{\ln(1/B)} \right\rceil}^+  \left\lfloor \frac{FM^2G^2}{2\bar\varepsilon^2\left[ B-\tfrac{1}{2F} \right]} \right\rfloor +  2  \left\lfloor 1 +\frac{\ln \left( \frac{1+\bar\varepsilon'}{B} \right)}{\ln(1/B)}  \right\rfloor \left\lfloor \frac{FM^2G^2}{2\varepsilon^2\left[ B-\tfrac{1}{2F} \right]} \right\rfloor  \\ \\

 \ \le &  2{\left[ \frac{\ln\left(1+\frac{f(x^0)-f^*}{f^* - \ua}\right)}{\ln(1/B)} \right]}  \left[ \frac{FM^2G^2}{2\bar\varepsilon^2\left[ B-\tfrac{1}{2F} \right]} \right] +  2  \left[ 2 +\frac{\ln \left(1+\bar\varepsilon' \right)}{\ln(1/B)}  \right] \left[ \frac{FM^2G^2}{2\varepsilon^2\left[ B-\tfrac{1}{2F} \right]} \right]  \\ \\

 \ \le & M^2 G^2 \left( 48.5 \ln\left( 1 + \tfrac{f(x^0)-f^*}{f^* - \ua} \right) \ +  \  18 \left(\frac{1+\varepsilon'}{\varepsilon'}\right)^2 \right) \ ,
\end{array}$$
where the second inequality follows from substituting in the values $B = 1/\sqrt{e}$, $F=\sqrt{e}$, and $\bar\varepsilon' = 0.9$, and rounding terms upward.  This last expression then is rounded upward to yield the desired iteration bound. \qed

\subsection{Subgradient Descent when $f^*$ is known}

In the special case when $f^*$ is known, we can obtain a computational guarantee that is of the same order as that of Theorem \ref{th1} by directly using the standard Subgradient Descent Method (Algorithm \ref{subgradientmethod}) with the (standard) step-size rule $ \alpha_i :=(f(x^i) - f^*)/\|g_i\|^2 $.  This is shown in the following theorem.\medskip

\begin{thm}{\bf (Complexity Bound for standard Subgradient Descent when $f^*$ is known)}\label{th3} Let the step-sizes for the Subgradient Descent Method (Algorithm \ref{subgradientmethod}) applied to solve problem \eqref{poi1} be chosen as:
$$ \alpha_i := \frac{f(x^i) - f^*}{\|g_i\|^2} \ , $$
and suppose that $N \ge 0$ and satisfies $$N \ \ge \ 2 M^2 G^2 \left(1 + 2.9  \ln\left(\frac{f(x^0)-f^*}{f^* - \ua} \right) + 2.9 \ln\left(\frac{1}{\varepsilon'}\right)  + 6.8\left(\frac{1}{\varepsilon'}\right) + 2\left(\frac{1}{\varepsilon'}\right)^2 \right) \ . $$
Then:
\begin{equation}\label{eq:known}
\ \ \ \ \ \ \\ \ \ \ \ \ \ \frac{f(x^N_b) - f^*}{f^*-\ua}\ \ \leq \ \varepsilon' \ . \ \ \ \ \ \ \ \ \ \ \ \ \ \ \ \qed
\end{equation}\end{thm}\medskip

The computational guarantee above is an almost-exact generalization of Theorem 3.7 of Renegar \cite{renegar2015}, which therein pertains to a specific transformed conic optimization problem.  The proof of this theorem follows the logic for the proof of Theorem 3.7 of \cite{renegar2015} in many respects as well.

Notice that up to an absolute constant, the computational guarantee of Theorem \ref{th3} is essentially the same as that of Theorem \ref{th1} in the worst case.\medskip

{\bf Proof of Theorem \ref{th3}:} We will presume that $\frac{f(x^0)-f^*}{f^*-\ua} >  \varepsilon'$, since otherwise \eqref{eq:known} is satisfied trivially for all $N \ge 0$.  Let $B \in (0,1)$ be a given fractional quantity.  Define $K_0:=0$, and for all $i$ such that $f^{K_{i}}_b - f^* > 0$ define $K_{i+1}$ inductively as the smallest iteration index of Subgradient Descent for which $f^{K_{i+1}}_b - f^* \le B(f^{K_{i}}_b - f^*)$.  Notice that so long as $f^{K_{i}}_b - f^* > 0$ it follows using part $(ii)$ of Theorem \ref{subgrad} that $K_{i+1}$ exists (i.e., is finite).  Let $i'$ be the smallest sub-index $i$ for which $\frac{f^{K_{i'}}_b-f^*}{f^*-\ua} \le  \varepsilon'$.  It follows from the initial presumption above that $i'\ge 1$, and it holds for any $i \ge 0$ satisfying $i < i'$ that $\varepsilon'(f^* - \ua)  < f(x^{K_{i}})-f^* \le B^{i}(f(x^{K_0})-f^*)= B^{i}(f(x^0)-f^*)$, from which it follows that $i$, and hence also $i'$, is finite.  Furthermore, it holds for any $i \ge 0$ satisfying $i < i'$ that:
\begin{equation}\label{sat2}
\varepsilon'(f^* - \ua)  < f(x^{K_{i'-1}})-f^* \le B^{i'-1-i}(f(x^{K_i})-f^*) \ .
\end{equation} Using $i=0$ in \eqref{sat2} and taking logarithms yields:
\begin{equation}\label{sat3}
i'  < 1 + \frac{\ln\left(\frac{1}{\varepsilon'}\right) + \ln\left(\frac{f(x^0)-f^*}{f^* - \ua}\right)}{\ln(B^{-1})} \ .\end{equation} 
If $K_{i+1}$ exists (i.e., is finite), then it follows from part $(ii)$ of Theorem \ref{subgrad} that:
\begin{equation*}
f(x_b^{K_{i+1}-1})-f^* \le \frac{M\dist(x^{K_i}, \optim)}{\sqrt{K_{i+1} -1-K_i +1}} \ .
\end{equation*}
This last inequality can be rearranged to yield:
\begin{equation}\label{sat1}
K_{i+1} - K_i \le \displaystyle\frac{M^2\dist(x^{K_i}, \optim)^2}{\left(f(x_b^{K_{i+1}-1})-f^*\right)^2}
< \displaystyle\frac{B^{-2}M^2G^2\left(f(x^{K_{i}}) - \ua\right)^2}{\left(f(x^{K_{i}})-f^*\right)^2}
= B^{-2}M^2G^2\left(1+ \displaystyle\frac{f^* - \ua}{f(x^{K_{i}})-f^*}\right)^2
\end{equation} where the second inequality uses the definition of the growth constant $G$ as well as the fact that $f(x^{K_{i+1}-1}_b) - f^* > B(f(x^{K_{i}}) - f^*)$. Now putting all of this together we obtain: \begingroup
\addtolength{\jot}{1.1em}\begin{align*}
K_{i'} &= \  \displaystyle\sum_{i=0}^{i'-1}\left( K_{i+1} - K_i \right)\\
&\le \   B^{-2}M^2G^2 \displaystyle\sum_{i=0}^{i'-1} \left(1+ \displaystyle\frac{f^* - \ua}{f(x^{K_{i}})-f^*}\right)^2\\
&\le \  B^{-2}M^2G^2 \displaystyle\sum_{i=0}^{i'-1} \left(1+ \displaystyle\frac{1}{\varepsilon'}B^{i'-1-i}\right)^2  \\
& = \  B^{-2}M^2G^2 \displaystyle\sum_{j=0}^{i'-1} \left(1+ \left(\frac{2}{\varepsilon'}\right)B^{j}+ \left(\frac{1}{\varepsilon'}\right)^2(B^2)^{j}\right)\\
& \le \   B^{-2}M^2G^2  \left(i'+ \left(\displaystyle\frac{2}{\varepsilon'}\right)\displaystyle\frac{1}{1-B}+ \left(\displaystyle\frac{1}{\varepsilon'}\right)^2\displaystyle\frac{1}{1-B^2}\right)  \\
& \le\   B^{-2}M^2G^2  \left(1 + \displaystyle\frac{\ln\left(\displaystyle\frac{1}{\varepsilon'}\right) + \ln\left(\displaystyle\frac{f(x^0)-f^*}{f^* - \ua}\right)}{\ln(B^{-1})} + \left(\displaystyle\frac{2}{\varepsilon'}\right)\displaystyle\frac{1}{1-B}+ \left(\displaystyle\frac{1}{\varepsilon'}\right)^2\displaystyle\frac{1}{1-B^2}\right) \ ,
\end{align*}\endgroup  where the first inequality is from \eqref{sat1}, the second inequality uses \eqref{sat2}, the third inequality replaces the two finite geometric series with corresponding infinite series, and the fourth inequality uses \eqref{sat3}.  Finally, using the value of $B = 1/\sqrt{2}$ and substituting into the above yields the result.\qed


We remark that one obtains the precise constants of Theorem 3.7 of \cite{renegar2015} by using $B=1/2$.  Choosing $B$ to optimize the absolute constant of the $(1/\varepsilon')^2$ term yields $B=1/\sqrt{2}$ and the absolute constants as presented in the statement of the threorem.  Choosing $B$ to optimize the absolute constant of the $\ln\left(\frac{f(x^0)-f^*}{f^* - \ua}\right)$ term would yield $B=1/\sqrt{e}$ with the coefficient of $2$ in the $\ln(\cdot)$ terms.

\subsection{Non-Smooth Optimization using a New Smooth Approximations Method}\label{generalissimo}

As first proposed by Nesterov \cite{nest05smoothing}, there are many practical settings wherein one can approximate the non-smooth convex function $f(\cdot)$ by a smooth convex function $f_\mu(\cdot)$, where the sense of the approximation depends on the parameter $\mu$.  If the smooth approximation $f_\mu(\cdot)$ is computationally easy to work with, one can then use the Accelerated Gradient Method (Algorithm \ref{accsimplegradientscheme}) to approximately optimize $f_\mu(\cdot)$ (thereby also approximately optimizing $f(\cdot)$) on the feasible set $Q$.  There are a variety of techniques that can be used to construct a parametric family of smooth functions $f_\mu(\cdot)$ depending on the known structure of $f(\cdot)$ and $Q$, see \cite{nest05smoothing} as well as \cite{nest07smoothsdp} and Beck and Teboulle \cite{btsmoothing} among others.  For our purposes herein, we will suppose that there is a smoothing technique with the following two properties:\medskip

\noindent (i) there is a known constant $\bar D > 0$ such that for any given $\mu >0$ we can construct a smooth convex function $f_\mu(\cdot):  Q \rightarrow \mathbb{R}$ which is not far from $f(\cdot)$, namely:
\begin{equation}\label{cliffman}
f(x) - \bar D \mu \ \le \ f_\mu(x)  \  \le \  f(x) \  \ \  \mbox{for~all~} x \in Q \ , \ \mathrm{and}
\end{equation}
\noindent (ii) $f_\mu(\cdot)$ has Lipschitz continuous gradient on $Q$ with Lipschitz constant $L_\mu$ satisfying
\begin{equation}\label{cliffman2} L_\mu \le \bar A/\mu\end{equation} for some known positive constant $\bar A$.

These properties can be used to design an implementation of the Accelerated Gradient Method (Algorithm \ref{accsimplegradientscheme}) applied to $f_\mu(\cdot)$, that can be used to compute an absolute $\varepsilon$-optimal solution of the original optimization problem \eqref{poi1}.  Using the scheme developed in \cite{nest05smoothing} in conjunction with the Accelerated Gradient Method (Algorithm \ref{accsimplegradientscheme}) yields an iteration complexity bound of
\begin{equation}\label{papa} \check N := \left\lceil \frac{\sqrt{8\bar A\bar D}\dist(x^0, \optim)}{\varepsilon} -1 \right\rceil
\end{equation} to obtain an (absolute) $\varepsilon$-optimal solution of \eqref{poi1} for a suitably designed version of the basic method \cite{nest05smoothing}.

Herein we develop a variant of the basic smoothing method to solve the optimization problem \eqref{poi1} that yields a new computational guarantee that can improve on \eqref{papa} in many cases.  Algorithm \ref{smoothedscheme} presents parametric smoothing and restarting method for computing an $\varepsilon'$-relative solution of the optimization problem \eqref{poi1} for the non-smooth objective function $f(\cdot)$ based on successive smooth approximations and re-starting of the Accelerated Gradient Method (Algorithm \ref{accsimplegradientscheme}).  In the description of Algorithm \ref{smoothedscheme} the general notation ``$x_{i,j} \gets \mathrm{AGM}(f_{\mu}(\cdot), \ x_{i,0}, \ j)$'' denotes assigning to $x_{i,j}$ the $j^{\mathrm{th}}$ iterate of the Accelerated Gradient Method applied to the optimization problem \eqref{poi1} with objective function $f_{\mu}(\cdot)$ using the initial point $x_{i,0} \in Q$.\medskip

\begin{algorithm}
\caption{Parametric Smoothing/Restarting Method using $f_\mu(\cdot)$}\label{smoothedscheme}
$ \ $
\begin{algorithmic}
\STATE {\bf Initialize.}  Initialize with $x^0 \in Q$ and $\varepsilon' >0$ .  \\
Define $ B := \tfrac{1}{2}$ , $t := \tfrac{1}{8}$ . \\
Set $x_{1,0} \gets x^0$ , $i \gets 1$ .
$ \ $ \\
$ \ $ \\
At outer iteration $i$ :
\STATE 1. {\bf Set smoothing parameters.} $\mu^1_i \gets \displaystyle\frac{t \cdot (f(x_{i,0}) - \ua)}{\bar D}$ , $\mu^2_i \gets \displaystyle\frac{t \varepsilon' \cdot (f(x_{i,0}) - \ua)}{\bar D}$ .\\
\STATE 2. {\bf Initialize inner iteration.} $K_i \gets +\infty$ , $j \gets 0$  \\\medskip
\STATE 3. {\bf Run inner iterations.} At inner iteration $j$: \\
\ \ \ \ \ \ \ \ \ \ \ \ \ \ \ \ \ \ \ \ \ \ \ \ \ \ \ (3a.) If \ \ $\displaystyle\frac{f(x_{i,j})-\ua }{f(x_{i,0})-\ua} \ge B$ , then \\\medskip
\ \ \ \ \ \ \ \ \ \ \ \ \ \ \ \ \ \ \ \ \ \ \ \ \ \ \ \ \ \ \ \ \ \ \ \ \ \ \ \ \ \ $x_{i,j+1} \gets \mathrm{AGM}(f_{\mu^1_i}(\cdot), \ x_{i,0}, \ j+1)$ , \\\medskip
\ \ \ \ \ \ \ \ \ \ \ \ \ \ \ \ \ \ \ \ \ \ \ \ \ \ \ \ \ \ \ \ \ \ \ \ \ \ \ \ \ \ $y_{i,j+1} \gets \mathrm{AGM}(f_{\mu^2_i}(\cdot), \ x_{i,0}, \ j+1)$ , \\\medskip
\ \ \ \ \ \ \ \ \ \ \ \ \ \ \ \ \ \ \ \ \ \ \ \ \ \ \ \ \ \ \ \ \ \ \ \ \ \ \ \ \ \ $j \gets j+1$, and Goto (3a.) \\\medskip
\ \ \ \ \ \ \ \ \ \ \ \ \ \ \ \ \ \ \ \ \ \ \ \ \ \ \ \ \ \ \ \ \ \ \ Else $K_i \gets j$, $x_{i+1,0} \gets x_{i,j}$, $i \gets i+1$, and Goto Step 1. \\\medskip
\end{algorithmic}
\end{algorithm}

At the $i^{\mathrm{th}}$ outer iteration of Algorithm \ref{smoothedscheme}, the algorithm sets two different smoothing parameters in Step (1.), namely $\mu^1_i$ and $\mu^2_i$, where $\mu_i^2$ differs from $\mu_i^1$ by the relative accuracy input value $\varepsilon'$.  The algorithm then runs the Accelerated Gradient Method with starting point $x_{i,0}$ simultaneously on the two smoothed functions $f_{\mu_i^1}(\cdot)$ and $f_{\mu_i^2}(\cdot)$, using the double indexing notation of $x_{i,j}$ and $y_{i,j}$ to denote iteration $j$ of the Accelerated Gradient Method initialized at the point $x_{i,0}$ for optimizing $f_{\mu_i^1}(\cdot)$ and $f_{\mu_i^2}(\cdot)$ on $Q$, respectively.  Notice that the smoothing parameters $\mu_i^1$ and $\mu_i^2$ decrease over the course of the outer iterations, as it makes more sense to set these values higher at first and then decrease them as the solution is approached.  The outer iteration $i$ runs until the ratio test in Step (3a.) fails, at which point the current point $x_{i,j}$ becomes the starting point of the next outer iteration, namely $x_{i+1,0} \gets x_{i,j}$.  The counter $K_i$ records the number of inner iterations $j$ of outer iteration $i$.  Regarding counting of iterates computed in Algorithm \ref{smoothedscheme}, we will say that the algorithm has computed an iterate whenever it calls $\mathrm{AGM}( \cdot, \cdot, \cdot)$.  There are therefore two computed iterates at each inner iteration.  

Restarting for accelerated gradient methods for strongly convex functions has been studied in \cite{o2013adaptive} and \cite{su2014differential}.  To the best of our knowledge, restarting of accelerated methods in the absence of strong convexity was first used in Renegar \cite{renegar2014}, and Algorithm \ref{smoothedscheme} exploits this and other ideas from \cite{renegar2014} and \cite{renegar2016b} as well.  We have the following computational guarantee associated with Algorithm \ref{smoothedscheme}. \medskip

\begin{thm}{\bf (Complexity Bound for Parametric Smoothing/Restarting Method (Algorithm \ref{smoothedscheme}) for Non-smooth Optimization)}\label{th4} Suppose that $f_\mu(\cdot)$ satisfies the smoothing conditions \eqref{cliffman} and \eqref{cliffman2}.  Within a total number of computed iterates that does not exceed $$23 G\sqrt{\bar A\bar D}\left(1 + 1.42\ln\left( 1+ \frac{f(x^0) - f^* }{f^* - \ua }\right)  + 2\left(\frac{1}{\varepsilon'} \right)\right) \ , $$
Algorithm \ref{smoothedscheme} will compute an iterate $y_{i,j}$ for which \begin{equation*}
\frac{f(y_{i,j}) - f^*}{f^* - \ua} \le \varepsilon'  \  .   \qed \end{equation*} 
\end{thm}

Similar to Theorem \ref{th1}, the dependence in Theorem \ref{th4} on the quality of the initial iterate is logarithmic in the initial optimality gap $f(x^0)-f^*$.  Also, the factor involving $1/\varepsilon'$ in Theorem
\ref{th4} is independent of the quality of the initial iterate, unlike that of the standard bound for the smoothing method given in \eqref{papa}.  We will prove Theorem \ref{th4} by first establishing several propositions.  Throughout, for notational convenience, we will work with two constants $B$ and $t$ that must be chosen to satisfy
\begin{equation*}\label{pats} B > 0 \ , \ t > 0 \ , \ B - B^2 \ge 2t \ , \  \mbox{and} \ B \ge 4t \ , \end{equation*}

and whose specific values are set to $B = 1/2$ and $t=1/8$ in Algorithm \ref{smoothedscheme}.

The following proposition applies to the generic setting of the Accelerated Gradient Method applied to the smoothed function $f_\mu(\cdot)$.  Recall that $L_\mu$ denotes the Lipschitz constant of the gradient of $f_\mu(\cdot)$ on $Q$. \medskip

\begin{prop}{\label{prop_1}}  Given the smoothing parameter $\mu >0$ and a given constant $\beta>0$, define $Y:=\left\lceil G\sqrt{2\beta} -1 \right\rceil $.  Let $x_{k} \gets \mathrm{AGM}(f_\mu(\cdot), \hat x^0, k)$ denote the $k^{\mbox{th}}$ iterate of the Accelerated Gradient Method applied to the function $f_{\mu}(\cdot)$ with starting point $\hat x^0$.  For $k \geq Y$ it holds that:
\begin{equation}
f(x_{k}) - f^* \leq \frac{L_\mu}{\beta}(f(\hat x^0) - \ua)^2+ \mu \bar D \ .
\end{equation}
\end{prop}

{\bf Proof:} Note that for any $x\in \optim$ it holds that $f_{\mu}(x) \leq f^*$, whereby $\optim \subset S:= \{ x \in Q : f_{\mu}(x) \le f^* \}$.  It then follows from Theorem \ref{accgrad} applied to the function $f_{\mu}(\cdot)$ and using $\delta = f^* $ that for any $k\geq Y$ we have:
\begin{equation*}
\begin{array}{lcl}
f_{\mu}(x_{k}) - f^*  & \le &  \frac{2L_\mu}{(Y+1)^2}\dist(\hat x^0,S)^2\\ \\
& \le &  \frac{2L_\mu}{(Y+1)^2}\dist(\hat x^0,\optim)^2
 \le  \frac{2L_\mu}{(Y+1)^2}G^2\left(f(\hat x^0) - \ua \right)^2
 \le  \frac{L_\mu}{\beta}\left(f(\hat x^0) - \ua \right)^2 \ ,

\end{array}
\end{equation*}

%
where the second inequality uses the fact that $\optim \subset S$, the third inequality uses the definition of $G$, and the last inequality uses the value of $Y$.

Note from \eqref{cliffman} that $f(x) \le f_\mu(x) + \mu \bar D$, whereby:
\begin{equation*}
f(x_{k}) - f^* \le f_{\mu}(x_{k}) - f^* + \mu \bar{D} \le \frac{L_\mu}{\beta}\left(f(\hat x^0) - \ua \right)^2 + \mu \bar D  \ . \ \qed
\end{equation*}

We now apply Proposition \ref{prop_1} to the setting of the Parametric Smoothing/Restarting Method (Algorithm \ref{smoothedscheme}). \medskip

\begin{prop}{\label{prop_2}}  Let $i$ be the index of an outer iteration of Algorithm \ref{smoothedscheme}.
Define $T:=\left\lceil \frac{G\sqrt{2\bar A\bar{D}}}{t} -1 \right\rceil $.  If $k \ge T$ and $x_{i,k}$ exists, then it holds that:
\begin{equation*}
f(x_{i,k}) - f^* \leq 2 t (f(x_{i,0}) - \ua) \ .
\end{equation*}
\end{prop} 
{\bf Proof:}  The proof follows by applying Proposition \ref{prop_1} with $\mu = \mu^1_i = \frac{t \cdot (f(x_{i,0}) - \ua)}{\bar D}$, $\beta =\frac{\bar A \bar D}{t^2}$, $Y=T$, and $\hat x^0 = x_{i,0}$.  It then follows that
$$f(x_{i,k}) - f^* \leq \frac{L_\mu}{\beta}\left(f(\hat x^0) - \ua \right)^2 + \mu \bar D \le \frac{\bar A}{\mu \beta}\left(f(\hat x^0) - \ua \right)^2 + \mu \bar D = 2 t (f(x_{i,0}) - \ua) \ ,
$$ where the second inequality uses $L_\mu \le \bar A/\mu$ from \eqref{cliffman2} and from substituting in the values of $\mu$ and $\beta$.  \qed \medskip

\begin{prop}{\label{prop_3}} Let $i$ be the index of an outer iteration of Algorithm \ref{smoothedscheme}.  Define $U:=\left\lceil \frac{G\sqrt{2\bar A\bar{D}}}{\varepsilon' t} -1 \right\rceil $.  If $k \ge U$ and $y_{i,k}$ exists, then it holds that:
\begin{equation*}
f(y_{i,k}) - f^* \leq 2\varepsilon' t (f(x_{i,0}) - \ua) \ .
\end{equation*}
\end{prop}

{\bf Proof:}  The proof follows by applying Proposition \ref{prop_1} with $\mu = \mu^2_i = \frac{t \varepsilon' \cdot (f(x_{i,0}) - \ua)}{\bar D}$, $\beta =\frac{\bar A \bar D}{t^2(\varepsilon')^2}$, $Y=U$, and $\hat x^0 = x_{i,0}$.  It then follows that
$$ f(y_{i,k}) - f^* \leq \frac{L_\mu}{\beta}\left(f(\hat x^0) - \ua \right)^2 + \mu \bar D \le \frac{\bar A}{\mu \beta}\left(f(\hat x^0) - \ua \right)^2 + \mu \bar D = 2 t \varepsilon'(f(x_{i,0}) - \ua) \ ,
$$ where the second inequality uses $L_\mu \le \bar A/\mu$ from \eqref{cliffman2} and the final equality derives from substituting in the values of $\mu$ and $\beta$.  \qed \medskip

The next three propositions pertain to Algorithm \ref{smoothedscheme} as well as to a more general setting which will be used in Section \ref{mtc} to prove computational guarantees for algorithms when $f(\cdot)$ is smooth.  The more general setting is described in the body of the following proposition. \medskip

\begin{prop}\label{prop_4}
Let $B,v>0$ be constants satisfying $B-B^2 \ge v$, $B\ge 2v$. Consider an algorithm with outer and inner iterations indexed with counters $i$ and $j$, respectively (such as Algorithm \ref{smoothedscheme}), with initial iterate $x^0$ that is used to set $x_{1,0} = y_{1,0} \gets x^0$ in simultaneous running of the Accelerated Gradient Method using the same indexing notation as in Algorithm \ref{smoothedscheme}, and where $x_{i+1,0} = y_{i+1,0} \gets x_{i,K_i}$ where $K_i \gets j$ denotes the first index $j$ for which $\frac{f(x_{i,j})-\ua}{f(x_{i,0})-\ua}<B$. Suppose that there are sequences $\{J_i\}$ and $\{I_i\}$ indexed over the outer iteration counter $i$ such that the following conditions are satisfied:

\ \ \ \ \ \ (i) for all $k\ge J_i$ it holds that $f(x_{i,k}) - f^* \le v\left(f(x_{i,0}) - \ua \right)$, and

\ \ \ \ \ \ (ii) for all $k\ge I_i$ it holds that $f(y_{i,k}) - f^* \le v \varepsilon' \left(f(x_{i,0}) - \ua \right)$.

Let $p$ denote the number of outer iterations $i$ for which $K_i$ is finite.  Then
\begin{equation*}
p \le \left\lfloor \frac{\ln \left( 1+ \frac{f(x^0) - f^* }{f^* - \ua }\right)}{\ln(1/B)} \right\rfloor \ .
\end{equation*} Furthermore, if $i \ge 1$ and $i \le p-1$, then $K_i \le J_i$.
\end{prop}

{\bf Proof:}  If $p=0$ the results follow trivially, so let us suppose that $p \ge 1$, whereby $K_p$ is finite and  $x_{p+1,0}$ exists.  It then follows that $f^* - \ua \le f(x_{p+1,0}) - \ua \le B^{p}(f(x_{1,0})-\ua)$, and taking logarithms yields the proof of the bound on $p$.

Suppose additionally that $i \ge 1$ and $i \le p-1$.  Let us assume that $K_i \ge J_i +1$, from which we will derive a contradiction.  We have $$  f(x_{i,K_i -1})-f^*  \le  v(f(x_{i,0}) - \ua)  \le  (B-B^2)(f(x_{i,0}) - \ua) \ , $$ where the first inequality uses condition $(i)$ and the second inequality uses $B-B^2 \ge v$.    Also, $i+2 \le p+1$, whereby $x_{i+2,0}$ exists and therefore satisfies $f^* - \ua \le f(x_{i+2,0}) - \ua \le B^2(f(x_{i,0}) - \ua) $.  Combining this inequality with that above yields $f(x_{i,K_i -1})-f^*  \le B(f(x_{i,0}) - \ua) - f^* + \ua $, which rearranges to yield:
$$\frac{f(x_{i,K_i -1})-\ua}{f(x_{i,0})-\ua} \le B \ , $$ and which contradicts the definition of $K_i$.  Therefore $K_i \le J_i$. \qed \medskip

\begin{prop}{\label{prop_5}} Under the same setting, notation, and conditions $(i)$ and $(ii)$ of Proposition \ref{prop_4}, let $i$ be the index of an outer iteration. If $j \ge J_i$ and $x_{i,j+1}$ exists, then:
$$\frac{f(x_{i,0})-\ua}{f^*-\ua} \le \frac{1}{B-v} \ .$$ Furthermore, if also $j \ge \max\{J_i,I_i\}$, then
$$\frac{f(y_{i,j}) - f^*}{f^* - \ua} \le \varepsilon' \ .
$$\end{prop}

{\bf Proof:}  Since $j \ge J_i$ it follows from condition $(i)$ that
\begin{equation*}
f(x_{i,j}) - f^* \le v( f(x_{i,0}) - \ua ) \ ,
\end{equation*}
and also since $x_{i,j+1}$ exists then $K_i \ge j+1$, whereby:
\begin{equation*}
\frac{f(x_{i,j})-\ua}{f(x_{i,0})-\ua} \ge B \ .
\end{equation*}
It then follows from these two inequalities that
\begin{equation*}\label{jill}
\frac{f(x_{i,0})-\ua}{f^*-\ua} = \frac{1}{\frac{f(x_{i,j})-\ua}{f(x_{i,0})-\ua} - \frac{f(x_{i,j})-f^*}{f(x_{i,0})-\ua}} \le \frac{1}{B-v} \ .
\end{equation*}
If also $j\ge I_i$, then we have from condition $(ii)$ that
$$ f(y_{i,j}) - f^* \le v \varepsilon' (f(x_{i,0})-\ua) \le \frac{v \varepsilon' }{B - v}(f^* - \ua)\le (f^* - \ua)\varepsilon' \ , $$ where the first inequality is from condition $(ii)$, the second inequality uses \eqref{jill}, and the third inequality uses $B \ge 2v$. \qed \medskip

\begin{prop}\label{prop_6}
Under the same setting, notation, and conditions $(i)$ and $(ii)$ of Proposition \ref{prop_4}, let $\hat N$ count the total number of inner iterations prior to and including the first iteration for which $y_{i,j}$ is an $\varepsilon'$-relative solution \eqref{asp}.  Then
\begin{equation}\label{kaylin} \mbox{either} \ \ \ (i) \  \hat N \le \sum_{i=1}^{p+1} J_i +  I_{p+1}\ , \ \ \ \ \ \ \mbox{or} \ \ \ \ \ \ (ii)\  \hat N \le \sum_{i=1}^{p+1} J_i + I_{p} + I_{p+1} \  \mbox{and} \  K_p \ge J_p+1\ .
\end{equation}
\end{prop}

{\bf Proof:}
First consider the case when $p=0$.  Then $K_1 = +\infty$ and therefore with $i=1$ we have $x_{i,j+1}$ exists for $j = \max\{J_1,I_1\}$, whereby from Proposition \ref{prop_5} it holds that $y_{1,j}$ satisfies \eqref{asp}.   In this case $\hat N \le  j =  \max\{J_1,I_1\} \le  J_1 +  I_1 =  \sum_{i=1}^{p+1} J_i +   I_{p+1} $ and therefore $(i)$ of \eqref{kaylin} is satisfied.

Next consider the case where $p\ge 1$ and $K_p \ge \max\{J_p, I_p\}+1$.  Let $i$ be the index of an outer iterate.  If $i \le p-1$ it follows from Proposition \ref{prop_4} that $K_i \le J_i$.  For $i=p$ it holds for this case that $K_p \ge \max\{J_p,I_p\}+1$, and it follows from Proposition \ref{prop_5} that $x_{p,j+1}$ exists for $j = \max\{J_p,I_p\}$ and therefore $y_{p,j}$ satisfies \eqref{asp}.  In this case $\hat N \le  \sum_{i=1}^{p-1}K_i +  \max\{J_p,I_p\} \le  \sum_{i=1}^{p-1}J_i +  \max\{J_p,I_p\} \le  \sum_{i=1}^{p}J_i +  I_p$ and $K_p\ge J_p + 1$ whereby $(ii)$ of \eqref{kaylin} is satisfied.

Next consider the case where $p\ge 1$ and $K_p \le \max\{J_p, I_p\}$ and also $K_p \le J_p$.  Let $i$ be the index of an outer iterate.  If $i \le p-1$ it follows from Proposition \ref{prop_4} that $K_i \le J_i$.  Since  $K_{p+1} = +\infty$ it follows that $x_{p+1,j+1}$ exists for $j = \max\{J_{p+1},I_{p+1}\}$, whereby from Proposition \ref{prop_5} we have $y_{p+1,j}$ satisfies \eqref{asp}.  And since $K_p\le J_p$ in this case, it follows that $\hat N \le  \sum_{i=1}^{p-1}K_i +  J_{p} +  \max\{J_{p+1},I_{p+1}\} \le  \sum_{i=1}^{p+1} J_i +   I_{p+1}$, and therefore $(i)$ of \eqref{kaylin} is satisfied.

The last case is where $p\ge 1$ and $K_p \le \max\{J_p, I_p\}$ and also $K_p \ge J_p +1$.  Then just as in the third case above, we arrive at $\hat N \le  \sum_{i=1}^{p-1}K_i +  \max\{J_{p},I_{p}\} +  \max\{J_{p+1},I_{p+1}\} \le  \sum_{i=1}^{p+1} J_i +  I_p +  I_{p+1}$, and thus $(ii)$ of \eqref{kaylin} is satisfied, thereby proving \eqref{kaylin}. \qed\medskip

{\bf Proof of Theorem \ref{th4}}:  Algorithm \ref{smoothedscheme} satisfies the setting of Proposition \ref{prop_4}, and it follows from Propositions \ref{prop_2} and \ref{prop_3} that Algorithm \ref{smoothedscheme} satisfies conditions $(i)$ and $(ii)$ of Proposition \ref{prop_4} by letting $v=2t$, $J_i=T$, and $I_i=U$ for all outer iterations $i$.   Therefore the conclusions of Propositions \ref{prop_4}, \ref{prop_5}, and \ref{prop_6} all hold true.  Let $N$ denote the total number of iterates of Algorithm \ref{smoothedscheme} computed prior to and including the first iterate $y_{i,j}$ that is an $\varepsilon'$-relative solution \eqref{asp}.  Since two iterates are computed at each iteration, we have $N=2\hat N$ (where $\hat N$ is defined in Proposition \ref{prop_6}) and it follows from Proposition \ref{prop_6} that
$N =2 \hat N \le 2 \sum_{i=1}^{p+1} J_i + 2 I_p + 2 I_{p+1}$, since the right-side of this inequality dominates both bounds $(i)$ and $(ii)$ of \eqref{kaylin}.  Substituting in the values of $T$ and $U$ and the bound on $p$ from Proposition \ref{prop_4} we obtain:
$$\begin{array}{rcl}
N \ \le \ 2(p+1)T + 4 U \ & \le&  \  2\left\lfloor 1 + \frac{\ln \left( 1+ \frac{f(x^0) - f^* }{f^* - \ua }\right)}{\ln(1/B)} \right\rfloor
\displaystyle{\left\lceil \frac{G\sqrt{2\bar A\bar{D}}}{t} -1 \right\rceil} + 4 \displaystyle{\left\lceil \frac{G\sqrt{2\bar A\bar{D}}}{\varepsilon' t} -1 \right\rceil}\\ \\

 \ &\le& G\sqrt{\bar A\bar D} \displaystyle{\left( 22.63 + 32.65\ln\left( 1 + \tfrac{f(x^0)-f^*}{f^* - \ua} \right) \ +  \  45.26 \left(\frac{1}{\varepsilon'}\right) \right)} \ ,
\end{array}$$ where the third inequality follows from substituting in the values $B=\tfrac{1}{2}$ and $t=\tfrac{1}{8}
$, which then rounds up to the desired bound in the theorem. \qed

\section{Computational Guarantees when $f(\cdot)$ is Smooth}\label{mtc}

In this section we study the computational complexity of solving \eqref{poi1} in the case when $f(\cdot)$ is convex and differentiable on $Q$.  We assume that $\nabla f(\cdot)$ is Lipschitz on $Q$ as defined in \eqref{lipgrad}.

Let us first consider directly applying the Accelerated Gradient Method (Algorithm \ref{accsimplegradientscheme}) to solve \eqref{poi1}, and let us apply Theorem \ref{accgrad}.   Let $\varepsilon' >0$ denote the relative accuracy, and note again that an $\varepsilon'$-relative solution of \eqref{poi1} corresponds to an absolute $\varepsilon$-solution for $\varepsilon := \varepsilon' \cdot (f^* - \ua)$.  Let $x^0 \in Q$ be the initial point.  It then follows from Theorem \ref{accgrad} using $\delta = f^*$ (whereby $S_\delta = \{x \in Q : f(x) \le f^*\} = \optim$) that if \begin{equation}\label{football} N \ \ge  \ \frac{\sqrt{2}\sqrt{L}\dist(x^0, \optim)}{\sqrt{\varepsilon'}\sqrt{f^*-\ua}} \  -1 \ , \end{equation}
then
$$\frac {f(x^N) - f^*}{f^*-\ua}\ \ \leq \ \varepsilon' \ .
$$\medskip

Herein we will derive a new computational guarantee for a version of the Accelerated Gradient Method that can improve on \eqref{football} in many cases.  Our new version of the Accelerated Gradient Method periodically restarts the method with an appropriate rule for deciding when to do the restarts, and is presented in Algorithm \ref{smoothschemerestart}.  At the $i^{\mathrm{th}}$ outer iteration of Algorithm \ref{smoothschemerestart}  the algorithm starts the Accelerated Gradient Method at the point $x_{i,0}$ for optimizing $f(\cdot)$ on $Q$. The outer iteration $i$ runs until the ratio test in Step (2a.) fails, at which point the current point $x_{i,j}$ becomes the starting point of the next outer iteration, namely $x_{i+1,0} \gets x_{i,j}$. The counter $K_i$ records the number of inner iterations computed in outer iteration $i$. Similar to the notation in Algorithm \ref{smoothedscheme}, the notation ``$x_{i,j} \gets \mathrm{AGM}(f(\cdot), x_{i,0}, j)$'' in Algorithm \ref{smoothschemerestart} denotes assigning to $x_{i,j}$ the $j^{\mathrm{th}}$ iterate of the Accelerated Gradient Method applied
to the optimization problem \eqref{poi1} with objective function $f(\cdot)$ using the initial point $x_{i,0} \in Q$.\medskip

\begin{algorithm}
\caption{Accelerated Gradient Method with Simple Restarting}\label{smoothschemerestart}
$ \ $
\begin{algorithmic}
\STATE {\bf Initialize.}  Initialize with $x^0 \in Q$ .  \\
Define $B:=0.5$ \\
Set $x_{1,0} \gets x^0$ , $ i \gets 1$ .
$ \ $ \\
$ \ $ \\
At outer iteration $i$ :

\STATE 1. {\bf Initialize inner iteration.} $K_i\leftarrow + \infty$, $j \gets 0$  \\\medskip
\STATE 2. {\bf Run inner iterations.} At inner iteration $j$: \\
\ \ \ \ \ \ \ \ \ \ \ \ \ \ \ \ \ \  \ (2a.) If \ \ $\displaystyle\frac{f(x_{i,j})-\ua }{f(x_{i,0})-\ua} \ge B$ , then \\\medskip
\ \ \ \ \ \ \ \ \ \ \ \ \ \ \ \ \ \ \ \ \ \ \ \ \ \ \ \ \ \ \ \ \ \ \ \ \ \ \ \ \ \ $x_{i,j+1} \gets \mathrm{AGM}(f(\cdot), \ x_{i,0}, \ j+1)$ , \\\medskip
\ \ \ \ \ \ \ \ \ \ \ \ \ \ \ \ \ \ \ \ \ \ \ \ \ \ \ \ \ \ \ \ \ \ \ \ \ \ \ \ \ \ $j \gets j+1$, and Goto (2a.). \\\medskip
\ \ \ \ \ \ \ \ \ \ \ \ \ \ \ \ \ \ \ \ \ \ \ \ \ \ \ Else $K_i \gets j$, $x_{i+1,0} \gets x_{i,j}$, $i \gets i+1$, and Goto step 1. \\\medskip
\end{algorithmic}
\end{algorithm}\medskip

We have the following computational guarantee associated with Algorithm \ref{smoothschemerestart}. \medskip

\begin{thm}{\bf (Complexity Bound for Accelerated Gradient Method with Simple Restarting)}\label{th55} Within a total number of computed iterates that does not exceed $$
G\sqrt{L}\left(10\sqrt{f(x^0)-\ua}  \ + \ 12\left[\frac{\sqrt{f^*-\ua}}{\sqrt{\varepsilon'}}\right] \right) \ , $$
the Accelerated Gradient Method with Simple Restarting (Algorithm \ref{smoothschemerestart})  will compute an iterate $x_{i,j}$ for which \begin{equation*}
 \ \ \ \ \ \ \ \ \ \ \ \ \ \ \frac {f(x_{i,j}) - f^*}{f^*-\ua}\ \ \leq \ \varepsilon' \ . \ \ \ \ \ \ \ \ \ \ \ \ \ \ \ \qed
\end{equation*}\end{thm}

The computational guarantee in Theorem \ref{th55} can itself be bounded by:
\begin{equation}\label{leaves}
G\sqrt{L}\left(10\sqrt{f^*-\ua}  \ + 10\sqrt{\tfrac{L}{2}}\dist(x^0, \optim) \ + \ 12\left[\frac{\sqrt{f^*-\ua}}{\sqrt{\varepsilon'}}\right] \right) \ , \end{equation}
which follows from the chain of inequalities:
$$\begin{array}{rcl}\sqrt{ f(x^0) - \ua} &=& \sqrt{(f^* - \ua) + (f(x^0) - f^*)} \\ \\ &\le& \sqrt{(f^* - \ua) + \tfrac{L}{2}\dist(x^0, \optim)^2 } \  \le \ \sqrt{(f^* - \ua)} + \sqrt{\tfrac{L}{2}}\dist(x^0, \optim) \ . 
\end{array}$$ 
Comparing \eqref{leaves} with the standard bound for the Accelerated Gradient Method given in \eqref{football}, we see that the factor involving $1/\sqrt{\varepsilon'}$ in \eqref{leaves} is independent of $\dist(x^0, \optim)$, unlike the standard bound \eqref{football}.  Towards the proof of Theorem \ref{th55}, for notational convenience we will work with two constants $B$ and $v$ that must be chosen to satisfy
\begin{equation}\label{bucs} B > 0 \ , \ v > 0 \ , \ B - B^2 \ge v \ , \  \mbox{and} \ B \ge 2v \ , \end{equation}

and whose specific values are set to $B = 0.5$ in Algorithm \ref{smoothedscheme}, and $v=0.25$ . \medskip

\begin{prop}{\label{prop_222}}
Let $i$ be the index of an outer iteration of Algorithm \ref{smoothschemerestart}. Define \\ $J_i:=\left\lceil G\sqrt{\frac{2L(f(x_{i,0})-\ua)}{v}} -1 \right\rceil $. If $k \ge J_i$ and $x_{i,k}$ exists, then it holds that:
\begin{equation}
f(x_{i,k}) - f^* \leq v(f(x_{i,0}) - \ua) \ .
\end{equation}
\end{prop}
{\bf Proof:} It follows from Theorem \ref{accgrad} applied to the function $f(\cdot)$ and using $\delta = f^*  $ that for any $k\geq J_i$ we have:
\begin{equation*}
f(x_{i,k}) - f^* \le  \frac{2L}{(J_i+1)^2}\dist(x_{i,0},\optim)^2 \le  \frac{2L}{(J_i+1)^2}G^2\left(f(x_{i,0}) - \ua \right)^2 \le  v\left(f(x_{i,0}) -\ua \right) \ ,
\end{equation*}
where the second inequality uses the definition of $G$, and the last inequality uses the value of $J_i$. \qed \medskip

\begin{prop}{\label{prop_333}}
Let $i$ be the index of an outer iteration of Algorithm \ref{smoothschemerestart}. Define \\ $I_i:=\left\lceil G\sqrt{\frac{2L(f(x_{i,0})-\ua)}{v \varepsilon'}} -1 \right\rceil $. If $k \ge I_i$ and $x_{i,k}$ exists, then it holds that:
\begin{equation*}
f(x_{i,k}) - f^* \leq v \varepsilon'  (f(x_{i,0}) - \ua) \ .
\end{equation*}
\end{prop}
{\bf Proof:} The proof follows using identical logic as in Proposition \ref{prop_222}. \qed\medskip

%
%
%

 \medskip

{\bf Proof of Theorem \ref{th55}}:  Even though Algorithm \ref{smoothschemerestart} does not simultaneously run two versions of the Accelerated Gradient Method, we can still view Algorithm \ref{smoothschemerestart} as an instance of the general algorithm setting of Proposition \ref{prop_4} by simply defining $y_{i,j} := x_{i,j}$ for all $i,j$.  It follows from Propositions \ref{prop_222} and \ref{prop_333} that Algorithm \ref{smoothschemerestart} satisfies conditions $(i)$ and $(ii)$ of Proposition \ref{prop_4}, and therefore Propositions \ref{prop_4}, \ref{prop_5}, and \ref{prop_6} hold for Algorithm \ref{smoothschemerestart}.  Substituting in the values of $J_i$ and using the fact that $f(x_{i,0}) - \ua \le B^{i-1}(f(x_{1,0})-\ua)$ for all iteration counters $i$, we obtain:\begin{equation*}\label{cloud1}
\begin{array}{rcl}
\displaystyle\sum_{i=1}^{p+1} J_i& \le & \displaystyle\sum_{i=1}^{p+1} G\sqrt{\frac{2L(f(x_{i,0})-\ua)}{v}} \\ \\
 & \le & \left(\displaystyle\sum_{i=0}^{p} B^{\frac{i}{2}}\right)G\sqrt{\displaystyle\frac{2L(f(x_{1,0})-\ua)}{v}} \\ \\
& < & \left(\displaystyle\sum_{i=0}^{\infty} B^{\frac{i}{2}}\right)G\sqrt{\displaystyle\frac{2L(f(x_{1,0})-\ua)}{v}} \  =  \ \displaystyle\frac{G}{1-\sqrt{B}}\sqrt{\displaystyle\frac{2L(f(x^0)-\ua)}{v}} \ .
\end{array}
\end{equation*}

Next observe that $K_{p+1} = \infty \ge J_p$, whereby it follows from Proposition \ref{prop_5} with $i=p+1$ that $f(x_{p+1,0}) - \ua \le \frac{1}{B-v} (f^* - \ua)$, and therefore it holds that: $$I_{p+1} \le G\sqrt{\frac{2L(f(x_{p+1,0})-\ua)}{v\varepsilon'}} \le G\sqrt{\frac{2L(f^*-\ua)}{(B-v)v\varepsilon'}} \ . $$
Also, if $K_p \ge J_p+1$, then similarly applying Proposition \ref{prop_5} with $i=p$ using the logic above implies that $I_{p} \le  G\sqrt{\frac{2L(f^*-\ua)}{(B-v)v\varepsilon'}}$.

Let $N$ denote the total number of iterates of Algorithm \ref{smoothschemerestart} computed prior to and including the first iterate $x_{i,j}$ that is an $\varepsilon'$-relative solution \eqref{asp}.  Then $N = \hat N$ where $\hat N$ is defined in Proposition \ref{prop_6}.  In either case {\em (i)} or {\em (ii)} of \eqref{kaylin}, it follows from Proposition \ref{prop_6} that:
$$\begin{array}{rcl}
N \ = \ \hat N \ & \le& \displaystyle\sum_{i=1}^{p+1} J_i + I_{p+1} + G\sqrt{\frac{2L(f^*-\ua)}{(B-v)v\varepsilon'}} \\ \\

&\le & \displaystyle\frac{G}{1-\sqrt{B}}\sqrt{\frac{2L(f(x^0)-\ua)}{v}} + 2 G\sqrt{\frac{2L(f^*-\ua)}{(B-v)v\varepsilon'}}\\ \\

 \ &\le& G\sqrt{ L} \displaystyle{\left( 9.66 \sqrt{f(x^0)-\ua} \ +  \  11.32 \left(\frac{f^*-\ua}{\varepsilon'}\right) \right)} \ ,
\end{array}$$ where the third inequality follows from substituting in the values $B=\tfrac{1}{2}$ and $v=\tfrac{1}{4}$, which then rounds up to the bound stated in the theorem. \qed \medskip

It turns out that we can further improve the computational guarantee of Theorem \ref{th55} by further modifying the Accelerated Gradient Method with Simple Restarting (Algorithm \ref{smoothschemerestart}), if we know and can easily work with an adjoint representation of $f(\cdot)$ to do ``extra smoothing.''  Let us see how this can be done.  We will assume that $f(\cdot)$ has the representation:
\begin{equation}\label{adjoint}
f(x) = \displaystyle\max_{\lambda \in P} \{ \lambda^T Ax-d(\lambda)\} \ ,
\end{equation}
where $P$ is a convex set and $d(\cdot)$ is a strongly convex function on $P$ with strong convexity parameter $\sigma$ and for which $ \min_{\lambda \in P} d(\lambda) \ge 0$.  (See \cite{nesterovBook} for properties of strongly convex functions.)  It then follows that $f(\cdot)$ is a globally smooth convex function with Lipschitz constant at most $L :=\|A\|^2/\sigma$, see Nesterov \cite{nest05smoothing}.  We presume further that $A$, $d(\cdot)$, and $P$ are given and that the optimization problem in \eqref{adjoint} is simple to solve. That being the case, for a given $x \in Q$, if $\tilde \lambda$ solves the optimization problem \eqref{adjoint}, then it holds that $f(x) = \tilde \lambda^TAx - d(\tilde\lambda)$ and $\nabla f(x) = A^T\tilde\lambda$ .

In a similar spirit as the smoothing technique employed in Section \ref{generalissimo}, we will consider parametrically working with a modification $f_\mu(\cdot)$ of $f(\cdot)$ that is more smooth than $f(\cdot)$ by increasing the weight on the the strongly convex function $d(\cdot)$ in \eqref{adjoint}.  For any $\mu \ge 0$ define the function $f_\mu(\cdot)$ by:
\begin{equation}\label{adjointly}
f_\mu(x) = \displaystyle\max_{\lambda \in P} \{ \lambda^T Ax-(1+\mu)d(\lambda)\} \ .
\end{equation}
If $P$ is bounded, then $\bar D:=\max_{\lambda\in P} \{d(\lambda)\}$ is finite, and the above smoothing technique has the following two properties:

\noindent (i) $f_\mu(\cdot)$ is not far from $f(\cdot)$,
\begin{equation}\label{relations2}
f(x) - \bar D \mu \ \le \ f_\mu(x)  \  \le \  f(x) \  \ \  \mbox{for~all~} x \in Q \ , \ \mathrm{and}
\end{equation}
\noindent (ii) $f_\mu(\cdot)$ has Lipschitz continuous gradient on $Q$ with Lipschitz constant $L_\mu$ satisfying
\begin{equation}\label{cliffman4} L_\mu \le  L/(1+\mu) \ . \end{equation}

This setting is very similar to the properties we have for smoothing of a non-smooth function $f(\cdot)$ in Section \ref{generalissimo}, and the only difference is that the Lipschitz constant $L_{\mu}$ here is bounded above by $L/(1+\mu)$ instead of by $\bar{A}/\mu$ as was the case in \eqref{cliffman2}.


Let $\varepsilon' > 0$ be given.  As before, we aspire to compute an $\varepsilon'$-relative solution of \eqref{poi1} as defined in \eqref{asp}.  We will use and analyze the Parametric Smoothing/Rescaling Method (Algorithm \ref{smoothedscheme}) but with $f_\mu(\cdot)$ defined by \eqref{adjointly} and hence satisfying \eqref{relations2} and \eqref{cliffman4}.  We have the following computational guarantee associated with Algorithm \ref{smoothedscheme} applied to the case when $f(\cdot)$ is smooth and $f_\mu(\cdot)$ is given by \eqref{adjointly}. \medskip


\begin{thm}{\bf (Complexity Bound for Parametric Smoothing/Restarting Method (Algorithm \ref{smoothedscheme}) for Smooth Optimization)}\label{th5} Suppose that $f_\mu(\cdot)$ is given by \eqref{adjointly} and hence satisfies \eqref{relations2} and \eqref{cliffman4}.  Within a total number of computed iterates that does not exceed $$
G\sqrt{L\bar D}\left( 22.7 + 32.7 \ln \left( 1+ \frac{f(x^0) - f^* }{f^* - \ua }\right) + 32 \sqrt{\frac{f^*-\ua}{\varepsilon'}} \right) \ , $$
Algorithm \ref{smoothedscheme} will compute an iterate $y_{i,j}$ for which $$\frac{f(y_{i,j}) - f^*}{f^* - \ua} \le \varepsilon'  \  . \ \ \ \ \ \ \ \ \ \ \ \  \qed$$ 
\end{thm}

The dependence in Theorem \ref{th5} on the quality of the initial point is logarithmic in the optimality gap $f(x^0) - f^*$, while it is the square root of the optimality gap in Theorem \ref{th55}.  We will prove Theorem \ref{th5} by first proving two propositions.  For notational convenience we will work with two constants $B$ and $t$, whose specific values are $B = \tfrac{1}{2}$ and $t=\tfrac{1}{8}$.\medskip

\begin{prop}{\label{prop_22}}
Let $i$ be the index of an outer iteration of Algorithm \ref{smoothedscheme}. Define $T:=\left\lceil G\sqrt{\frac{2L\bar D}{t^2}} -1 \right\rceil $. If $k\ge T$ and $x_{i,k}$ exists, then:
\begin{equation*}
f(x_{i,k}) - f^* \leq 2t(f(x_{i,0}) - \ua) \ .
\end{equation*}
\end{prop}

{\bf Proof:} The proof follows by applying Proposition \ref{prop_1} with $\mu = \mu^1_i = \frac{t \cdot (f(x_{i,0}) - \ua)}{\bar D}$, $\beta =\frac{L \bar D}{t^2}$, $Y=T$, and $\hat x^0 = x_{i,0}$.  It then follows that
$$
f(x_{i,k}) - f^* \leq \frac{L_\mu}{\beta}\left(f(\hat x^0) - \ua \right)^2 + \mu \bar D \le \frac{L}{\mu \beta}\left(f( x_{i,0}) - \ua \right)^2 + \mu \bar D = 2 t (f(x_{i,0}) - \ua) \ ,
$$ where the second inequality uses $L_\mu \le L/(1+\mu) \le L/\mu$ from \eqref{cliffman4} and the final equality derives from substituting in the values of $\mu$ and $\beta$.  \qed \medskip

\begin{prop}{\label{prop_33}}
Let $i$ be the index of an outer iteration of Algorithm \ref{smoothedscheme}. Define \\ $I_i:=\left\lceil G\sqrt{\frac{2L(f(x_{i,0})-\ua)}{t\varepsilon' }} -1 \right\rceil $. If $k\ge I_i$ and $y_{i,k}$ exists, then:
\begin{equation*}
f(y_{i,k}) - f^* \leq 2t\varepsilon'(f(x_{i,0}) - \ua) \ .
\end{equation*}
\end{prop}
{\bf Proof:} The proof follows by applying Proposition \ref{prop_1} with $\mu = \mu^2_i = \frac{t \varepsilon' \cdot (f(x_{i,0}) - \ua)}{\bar D}$, $\beta =\frac{L (f(x_{i,0}) - \ua)}{t \varepsilon'}$, $Y=I_i$, and $\hat x^0 = x_{i,0}$.  It then follows that
$$
f(y_{i,k}) - f^* \leq \frac{L_\mu}{\beta}\left(f(\hat x^0) - \ua \right)^2 + \mu \bar D \le \frac{L}{ \beta}\left(f( x_{i,0}) - \ua \right)^2 + \mu \bar D = 2 t \varepsilon' (f(x_{i,0}) - \ua) \ ,
$$ where the second inequality uses $L_\mu \le L/(1+\mu) \le L$ from \eqref{cliffman4} and the final equality derives from substituting in the values of $\mu$ and $\beta$.  \qed \medskip

{\bf Proof of Theorem \ref{th5}}:  Algorithm \ref{smoothedscheme} satisfies the setting of Proposition \ref{prop_4}, and it follows from Propositions \ref{prop_22} and \ref{prop_33} that Algorithm \ref{smoothedscheme} satisfies conditions $(i)$ and $(ii)$ of Proposition \ref{prop_4} by letting $v=2t$ and $J_i=T$ for all outer iterations $i$.  Therefore the conclusions of Propositions \ref{prop_4}, \ref{prop_5}, and \ref{prop_6} all hold true.  Let $N$ denote the total number of iterates of Algorithm \ref{smoothedscheme} computed prior to and including the first iterate $y_{i,j}$ that is an $\varepsilon'$-relative solution \eqref{asp}.  Since two iterates are computed at each iteration, we have $N=2\hat N$, where $\hat N$ is defined in Proposition \ref{prop_6} and is bounded by either $(i)$ or $(ii)$ of \eqref{kaylin}.

Note that $K_{p+1} = \infty \ge T = J_i$, whereby it follows from Proposition \ref{prop_5} that $f(x_{p+1,0}) - \ua \le \frac{1}{B-2t} (f^* - \ua)$, and therefore $$I_{p+1} \le G\sqrt{\frac{2L(f(x_{p+1,0})-\ua)}{t\varepsilon'}} \le G\sqrt{\frac{2L(f^*-\ua)}{(B-2t)t\varepsilon'}}\  .$$
Similarly, if $K_p \ge T+1 = J_i + 1$, similar logic demonstrates that $I_{p} \le  G\sqrt{\frac{2L(f^*-\ua)}{(B-2t)t\varepsilon'}}$. Therefore, in either case $(i)$ or $(ii)$ of \eqref{kaylin} it holds that:

$$\begin{array}{rcl}
N \ = \ 2 \hat N \ & \le \ & 2\sum_{i=1}^{p+1} J_i + 2I_{p+1} + 2 G\sqrt{\frac{2L(f^*-\ua)}{(B-2t)t\varepsilon'}} \\ \\

\ & \le& 2(p+1)\left\lceil G\sqrt{\frac{2L\bar D}{t^2}} -1 \right\rceil + 2I_{p+1} + 2G\sqrt{\frac{2L(f^*-\ua)}{(B-2t)t\varepsilon'}} \\ \\

&\le & 2\left( 1 + \frac{\ln \left( 1+ \frac{f(x^0) - f^* }{f^* - \ua }\right)}{\ln(1/B)} \right)  G\sqrt{\frac{2L\bar D}{t^2}}  + 4 G\sqrt{\frac{2L(f^*-\ua)}{(B-2t)t\varepsilon'}}\\ \\

&\le & \ G\sqrt{L\bar D}\left( 22.7 + 32.7 \ln \left( 1+ \frac{f(x^0) - f^* }{f^* - \ua }\right) + 32 \sqrt{\frac{f^*-\ua}{\varepsilon'}} \right) \ ,
\end{array}$$ where the third inequality follows from substituting in the values $B=\tfrac{1}{2}$ and $t=\tfrac{1}{8}$, which then rounds up to the desired bound in the theorem. \qed
\appendix
\section{Appendix}

\subsection{Growth Constant $G$ and the Modulus of Weak Sharp Minima}\label{burke-ferris}

The optimal solution set $\optim$ of \eqref{poi1} is called a set of {\em weak sharp minima}
with modulus $\alpha$ if it holds that:
\begin{equation}
f(x)\geq f^*+\alpha\cdot \dist(x,\optim)  \ \ \mbox{for~all~} x \in Q \ . \label{eq:wsm}
\end{equation} This concept was first developed by Polyak \cite{polyak1979sharp} when $\optim$ is a singleton, and generalized by Burke and Ferris \cite{burke1993weak} to include the possibility of multiple optima.
The modulus of weak sharp minima has been a useful
tool in sensitivity analysis \cite{burke2001optimal,jourani2000hoffman}, convergence
analysis for certain problem classes \cite{burke1995gauss,burke1993weak}, linear regularity and error bounds \cite{burke2002weak,burke2005weak,burke2009weak}, perturbation properties of nonlinear optimization \cite{shapiro1988perturbation,shapiro1992perturbation,bonnans1994quadratic}, as well as in the finite termination of certain algorithms
\cite{polyakbook}, \cite{ferris1991finite}, and \cite{burke1995gauss}.

Comparing \eqref{eq:wsm} to \eqref{gg2}, we see that the modulus $\alpha$ of weak sharp minima is a close cousin of the growth constant $G$.  Indeed, if we were to loosen the restriction that $\ua$ be a strict lower bound and instead allow it to take the value $\ua = f^*$ in the definition of $G$ in \eqref{gg}, then we would obtain precisely that $G=\alpha^{-1}$.  However, the notion of $\ua$ being a strict lower bound is fundamental for the results herein.

Note that \eqref{eq:wsm} specifies the exact local growth of $f(\cdot)$ away from the set of optimal solutions.  And although as defined in \eqref{eq:wsm} the weak sharp minima is a global property, due to convexity it is essentially a local property and indeed its usefulness derives from the local nature of the weak sharp minima in a neighborhood of the optimal solution set.  This is in contrast to the growth constant $G$ as defined in \eqref{gg}, which by its nature is a global property as illustrated in the constructions in Figure \ref{fig1}.  Last of all we point out that while one can easily have $\alpha = 0$ for weak sharp minima (just let $f(x)$ be a differentiable convex function whose optimum is attained in the relative interior of $Q$), Theorem \ref{lu-theorem} shows that $G$ is finite for all reasonably-behaved convex functions.

\subsection{Proof of Theorem \ref{lu-theorem}}\label{app:lu-theorem}

{\bf Proof of Theorem \ref{lu-theorem}:}  Let us fix an optimal solution $x^* \in \optim$, and define $\delta:=\max_{v\in E_{\varepsilon}}  \| v-x^*\| $ and define $\bar{G} := \max\{\tfrac{\delta}{\varepsilon}, \tfrac{\delta}{f^*-\ua}\}$.  We will prove that for any $x \in Q$, the following inequality holds:
\begin{equation}\label{grow}
\dist\left(x,\optim\right)\leq \bar{G}\left(f(x)-\ua\right) \ ,
\end{equation}which then implies that $G \le \bar G$ is finite.  We consider two cases as follows:

\noindent {Case (i):  $x\in \optim_\varepsilon$.}  In this case we have $x = v+s$ where $v \in E_{\varepsilon}$ and $s \in S$.  Since $s$ is in the recession cone of $\optim_{\varepsilon}$ it holds that $x^* + s \in \optim$, whereby
\begin{equation}\label{church}
\dist(x, \optim)  \leq  \| x - (x^*+s) \| = \| v-x^*\| \le \delta \ ,
\end{equation} and therefore
\begin{equation*}
f(x) - \ua \geq f^* - \ua \ge \frac{(f^* - \ua)(\dist(x, \optim))}{\delta} \geq \bar G^{-1} \dist(x, \optim) \ ,
\end{equation*} which shows \eqref{grow} in this case.

\noindent {Case (ii):  $x \notin \optim_\varepsilon$.}  Let $x^1$ be the projection of $x$ onto $\optim$ and let $x^2$ be the point on the line segment from $x^1$ to $x$ that satisfies $f(x^2)= f^*+\varepsilon$.  (Existence of $x^2$ is guaranteed by continuity of $f(\cdot)$.)  Then
\begin{equation*}
f(x)-\ua \geq f(x) - f^* \geq  (f(x^2) - f^*)\frac{\| x-x^1 \|}{\| x^2 - x^1 \|} \geq  \frac{\varepsilon\| x - x^1\|}{\delta} = \frac{\varepsilon \dist(x, \optim)}{\delta} \ge \bar G^{-1} \dist(x, \optim)  \ ,
\end{equation*}

where the second inequality is from the convexity of $f(\cdot)$ which implies the chordal inequality $\frac{f(x) - f^*}{\| x-x^1 \|} \geq  \frac{f(x^2) - f^*}{\| x^2 - x^1 \|} $, and the third inequality uses $\| x^2-x^1 \| =\dist(x^2, \optim) \leq \delta$ (from \eqref{church}). The last equality above uses the fact that $\dist(x, \optim) = \| x-x^1 \|$.  This proves \eqref{grow} in this case. \qed\medskip


\bibliographystyle{amsplain}
\bibliography{GF-papers-nips-better}

\end{document}